\documentclass[11pt]{amsart}
\usepackage{amssymb,epic,eepic,verbatim,graphicx}
\newtheorem{theorem}{Theorem}[subsection]

\newtheorem{corollary}[theorem]{Corollary}

\newtheorem{proposition}[theorem]{Proposition}
\newtheorem{lemma}[theorem]{Lemma}
\newtheorem{lem}[theorem]{}
\theoremstyle{definition}
\newtheorem{definition}[theorem]{Definition}
\theoremstyle{remark}
\newtheorem{remark}[theorem]{Remark}
\newtheorem{example}[theorem]{Example}
\newcommand{\blem}{\begin{lem} \rm}
\newcommand{\elem}{\end{lem}}

\newcommand\A{\mathcal{A}}

\newcommand{\J}{\mathcal{J}}

\newcommand{\R}{\mathbb{R}}
\renewcommand{\H}{\mathbb{H}}

\newcommand{\C}{\mathbb{C}}

\newcommand{\Z}{\mathbb{Z}}

\newcommand{\ddd}{\frac{d}{d \delta}}

\newcommand\lie[1]{\mathfrak{#1}}

\newcommand{\g}{\lie{g}}

\newcommand{\on}{\operatorname}

\newcommand{\Def}{\on{Def}}

\newcommand{\pre}{{\on{pre}}}

\renewcommand{\vert}{{\on{vert}}}
\newcommand{\Edge}{\on{Edge}}

\newcommand{\End}{\on{End}}

\newcommand{\Aut}{ \on{Aut} } 
 
\newcommand{\aut}{ \on{aut} }

\renewcommand{\ker}{ \on{ker}}

\newcommand{\Vol}{  \on{Vol}}

\newcommand\dirac{/\kern-1.2ex\partial} 
\newcommand\qu{/\kern-.7ex/} 
\newcommand\lqu{\backslash \kern-.7ex \backslash} 

\newcommand\dr{r_+ \kern-.7ex - \kern-.7ex r_-}
 
\newcommand{\labell}\label

\renewcommand{\d}{{\on{d}}}
\newcommand{\ol}{\overline}
\newcommand{\olp}{\ol{\partial}}

\newcommand\eps{\epsilon}

\newcommand{\f}{\frac}

\newcommand{\hh}{{\f{1}{2}}}

\newcommand{\ti}{\tilde}

\newcommand\cE{\mathcal{E}}

\newcommand\cF{\mathcal{F}}
\newcommand\cI{\mathcal{I}}
\newcommand\cH{\mathcal{H}}

\newcommand\Map{\on{Map}}

\newcommand\ev{\on{ev}}

\newcommand\Vect{\on{Vect}}

\newcommand\G{\mathcal{G}}
\renewcommand\H{\mathcal{H}}
\newcommand\E{\mathcal{E}}
\newcommand\B{\mathcal{B}}
\renewcommand\Im{\on{Im}}

\newcommand\reg{{\on{reg}}}

\newcommand\bra[1]{ < \kern-.7ex {#1} \kern-.7ex >} 
\newcommand\bdefn{\begin{definition}}
\newcommand\edefn{\end{definition}}
\newcommand\bea{\begin{eqnarray*}}
\newcommand\eea{\end{eqnarray*}}
\newcommand\bcv{\left[ \begin{array}{r} }
\newcommand\ecv{\end{array} \right] }

\newcommand\bma{\left[ \begin{array}{l} }
\newcommand\ema{\end{array} \right]}
\newcommand\ben{\begin{enumerate}}
\newcommand\een{\end{enumerate}}
\newcommand\beq{\begin{equation}}
\newcommand\eeq{\end{equation}}
\newcommand\bex{\begin{example}}
\newcommand\bsj{\left\{ \begin{array}{rrr} }
\newcommand\esj{\end{array} \right\}}

\newcommand\mul{}

\newcommand\eex{\end{example}}

\newcommand\sx{*\kern-.5ex_X}

\def\mathunderaccent#1{\let\theaccent#1\mathpalette\putaccentunder}
\def\putaccentunder#1#2{\oalign{$#1#2$\crcr\hidewidth \vbox
to.2ex{\hbox{$#1\theaccent{}$}\vss}\hidewidth}}

\begin{document}

\title{Deformations of symplectic vortices}

\author{Eduardo Gonzalez} 

\address{
Department of Mathematics
University of Massachusetts Boston
100 William T. Morrissey Boulevard
Boston, MA 02125}
  \email{eduardo@math.umb.edu}

\author{Chris Woodward}

\address{Mathematics-Hill Center,
Rutgers University, 110 Frelinghuysen Road, Piscataway, NJ 08854-8019,
U.S.A.}  \email{ctw@math.rutgers.edu}

\thanks{Partially supported by NSF
 grant  DMS060509}

\begin{abstract}  
We prove a gluing theorem for a symplectic vortex on a compact complex
curve and a collection of holomorphic sphere bubbles.  Using the
theorem we show that the moduli space of regular stable symplectic
vortices on a fixed curve with varying markings has the structure of a
stratified-smooth topological orbifold.  In addition, we show that the
moduli space has a non-canonical $C^1$-orbifold structure.
\end{abstract}

\maketitle

\parskip 0in
\tableofcontents
\parskip .1in

\section{Introduction}

In this paper we generalize the following result on existence of
universal deformations for stable (pseudo-)holomorphic maps.  Let
$(X,\omega)$ be a compact symplectic manifold equipped with a
compatible almost complex structure $J$, and $(\mul{\Sigma},\mul{j})$
a compact nodal complex curve.  A map $\mul{u}: \mul{\Sigma} \to X$ is
{\em holomorphic} if
$$ \olp \mul{u}:= J_{\mul{u}} \circ \d \mul{u} - \d \mul{u} \circ \mul{j} = 0 $$
on each component of $\mul{\Sigma}$. One naturally has the notion of a
stratified-smooth {\em family} of holomorphic maps, and hence the
notion of a {\em deformation}, namely the germ of a family around the
central fiber together with an isomorphism of the central fiber with
the given map.  Recall that a deformation is {\em universal} if any
other deformation is obtained from it by pullback, in a unique way, by
a map of parameter spaces.  A holomorphic map $\mul{u}: \mul{\Sigma}
\to X$ is {\em regular} if the linearized Cauchy-Riemann operator is
surjective.  The following theorem is the result of the well-known
gluing construction for holomorphic maps, c.f. Ruan-Tian \cite{rt:hi}
or the text McDuff-Salamon \cite[Chapter 10]{ms:jh} in the case of
genus zero:

\begin{theorem}   \label{premain}  A regular holomorphic map $\mul{u}: \mul{\Sigma} \to X$ admits a stratified-smooth universal deformation iff it is stable.  
\end{theorem} 

The construction of the universal deformation proceeds via the
implicit function theorem.  For each element in the infinitesimal
deformation space of the stable map one first produces an approximate
solution and then applies the implicit function theorem to find an
exact solution.  Unfortunately one uses a different Sobolev space for
each ``gluing parameter'' controlling the domain, which means that it
is rather tricky to show that each nearby stable holomorphic map
occurs only once in the resulting family.  A slightly jazzed up
version of the above theorem implies that the gluing construction
gives rise to orbifold charts on the regular locus of the moduli space
of stable holomorphic maps.  Uniqueness of the universal deformations
implies that the smooth structures on each stratum are independent of
the Sobolev spaces used in the implicit function theorem.  One can
make these charts $C^1$-compatible by suitable choices of {\em gluing
  profiles}, that is, coordinates on the local deformation spaces;
however the $C^1$-structure on the moduli space is not canonical.  The
first part of the paper contains an exposition of the above theorem,
which is rather scattered in the literature.

The main result of the paper is a generalization of the theorem above
to certain {\em gauged (pseudo)holomorphic maps}, namely {\em
  symplectic vortices} as introduced by Mundet \cite{mun:ham} and
Cieliebak, Gaio and Salamon, see \cite{ciel:vor}.  Let $G$ be a
compact Lie group and $X$ a Hamiltonian $G$-manifold equipped with a
moment map $\Phi: X \to \g^*$ and an invariant almost complex
structure $J$.  Let $\Sigma$ be a compact smooth holomorphic curve
with complex structure $j$ and equipped with an area form
$\Vol_\Sigma$.  A {\em gauged holomorphic map} with values in $X$
consists of a smooth principal $G$-bundle $P \to \Sigma$, a connection
$A$ on $P$, and a smooth section $u: \Sigma \to P(X) := P \times_G X$
such that $\olp_A u = 0$ where $\olp_A$ is defined using the splitting
given by the connection $A$ and the complex structures $J, j$.  Let
$F_A \in \Omega^2(\Sigma, P(\g))$ denote the curvature of $A$ and
$P(\Phi): P(X) \to P(\g)$ the map induced by $\Phi$.  The space of
gauged holomorphic maps admits a formal symplectic structure depending
on a choice of invariant metric on $\g$ so that the action of the
group of gauge transformations is formally Hamiltonian.  A {\em
  symplectic vortex} is a pair in the zero level set of the moment
map: a pair $(A,u)$ such that
$$ \olp_A u = 0, \quad  F_A + u^* (P \times_G \Phi) \Vol_\Sigma = 0 .$$
Thus the moduli space $M(\Sigma,X)$ of symplectic vortices is the
symplectic quotient of the space of gauged maps by the group of gauge
transformations.  In certain cases where the moduli spaces are compact
Cieliebak, Gaio, Mundet, and Salamon \cite{ci:symvortex} and Mundet
\cite{mun:ham} constructed invariants that we will call {\em gauged
  Gromov-Witten invariants} by integration over these moduli spaces.
In general $M(\Sigma,X)$ admits a compactification $\ol{M}(\Sigma,X)$
consisting of {\em polystable symplectic vortices} given by allowing
$u$ to develop holomorphic sphere bubbles in the fibers of $P \times_G
X$.  A polystable vortex is {\em strongly stable} if the principal
component has finite automorphism group, and {\em regular} if a
certain linearized operator is surjective, that is, the moduli space
is formally smooth.  Our main result is the following:
\begin{theorem}  \label{main}  Let $\Sigma,X$ be as above.  
A regular strongly stable symplectic vortex from $\Sigma$ to $X$
admits a universal stratified-smooth deformation.
\end{theorem}  
\noindent 
Using the deformations constructed in Theorem \ref{main} we prove that
the moduli space $\ol{M}^{\reg}(\Sigma,X)$ of regular strongly stable
symplectic vortices admits the structure of an oriented
stratified-smooth topological orbifold, and (non-canonically) the
structure of a $C^1$-orbifold.  The first statement implies that if
$\ol{M}^\reg(\Sigma,X)$ is compact then it carries a rational
fundamental class.  The second statement implies for example, that if
the target carries a group action then the usual equivariant
localization theorems hold for the induced group action on the moduli
space.  In the case that $X$ is a smooth projective variety, algebraic
methods explain in \cite{cross} give similar results and provide
virtual fundamental classes on the moduli space.  However, the
symplectic gluing construction is interesting in its own right, not in
the least because it potentially extends to the case of Lagrangian
boundary conditions.  We understand that a forthcoming paper of Mundet
i Riera and Tian gives a gluing construction for two symplectic
vortices, when the structure group is the circle group.

\noindent {\em Acknowledgments:} We thank Ignasi Mundet i Riera, Melissa Liu,
and Robert Lipshitz for helpful comments and discussions.

\section{Deformations of holomorphic curves} 

The following section is essentially a review of the material that can
be found at the beginning of Siebert \cite{sie:gen}, with a few
additional comments incorporating terminology of Hofer, Wysocki, and
Zehnder \cite[Appendix]{ho:sc}.  In the first part we review the holomorphic
construction of universal deformations of stable curves. In the second
part, we study smooth deformations of curves.

\subsection{Holomorphic families of stable curves}  \label{stablecurves}

A compact, complex {\em nodal curve} $\mul{\Sigma}$ is obtained from a
collection $(\Sigma_1,\ldots,\Sigma_k)$ of smooth, compact, complex
curves by identifying a collection of distinct {\em nodal points}
$$\mul{w} = \{ \{ w_1^-,w_1^+ \}, \ldots, \{ w_m^-,w_m^+ \} \} .$$
For $l = 1,\ldots, m$, we denote by $\Sigma_{i^\pm(l)}$ the components
such that $w_l^\pm \in \Sigma_{i^\pm(l)}$.  A point $z \in
\mul{\Sigma}$ is {\em smooth} if it is not equal to any of the nodal
points.  A {\em marked nodal curve} is a nodal curve together with a
collection $\mul{z} = (z_1,\ldots, z_n)$ of distinct, smooth points.
An {\em isomorphism} of marked nodal curves $(\mul{\Sigma}_0,\mul{z}_0)$
to $(\mul{\Sigma}_1,\mul{z}_1)$ is an isomorphism $\phi: \mul{\Sigma}_0
\to \mul{\Sigma}_1$ of nodal curves such that $\phi(z_{0,i}) = z_{1,i}$
for $i = 1,\ldots, n$.  A marked nodal curve is {\em stable} if it has
finite automorphism group, that is, each component contains at least
three marked or nodal points if genus zero, or one special point if
genus one.  

The {\em combinatorial type} $\Gamma(\mul{\Sigma})$ of $\mul{\Sigma}$ is
the graph whose vertices are the components and edges are the nodes
and markings of $\mul{\Sigma}$.  The map $\mul{\Sigma} \mapsto
\Gamma(\mul{\Sigma})$ extends to a functor from the category of marked
nodal curves to the category of graphs. In particular, there is a
canonical homomorphism $\Aut(\mul{\Sigma}) \to
\Aut(\Gamma(\mul{\Sigma}))$, whose kernel is the product of the
automorphism groups of the components of $\mul{\Sigma}$.

Let $S$ be a complex variety (or scheme).  A {\em family of nodal
  curves} over $S$ is a complex variety $\mul{\Sigma}_S$ equipped with
a proper flat morphism $\pi: \mul{\Sigma}_S \to S$, such that each
fiber $\mul{\Sigma}_s, s \in S$ is a nodal curve.  A {\em deformation}
of a marked nodal curve $\mul{\Sigma}$ is a germ of a family of marked
nodal curves $\mul{\Sigma}_S$ over a pointed space $(S,0)$ together
with an isomorphism $\varphi: \mul{\Sigma}_0 \to \mul{\Sigma}$ of the
{\em central fiber} $\mul{\Sigma}_0$ with $\mul{\Sigma}$.  A
deformation $(\mul{\Sigma}_S, \varphi)$ of $\mul{\Sigma}$ is {\em
  versal} iff any other deformation $(\mul{\Sigma}'_S \to S',
\varphi')$ is induced from a map $\psi: S' \to S$ in the sense that
there exists an isomorphism $\phi$ of $\mul{\Sigma}'$ with the fiber
product $\mul{\Sigma_S} \times_S S'$ in a neighborhood of the central
fiber $ \mul{\Sigma}_0$.  A versal deformation is {\em universal} if
the map $\phi$ is the unique such map inducing the identity on
$\mul{\Sigma}_0$.  A deformation has {\em fixed type} if the
combinatorial type of the fiber is constant.  A {\em universal
  deformation of fixed type} is a deformation of fixed type, which is
universal in the above sense for deformations of fixed type.  The
space $\Def(\mul{\Sigma})$ of {\em infinitesimal deformations} of
$\mul{\Sigma}$ is the tangent space $T_0 S $ of the base $S$ of a
universal deformation, well-defined up to isomorphism.  We write
$\Def_\Gamma(\mul{\Sigma})$ for the space of infinitesimal
deformations of fixed type.  Let $\ti{\Sigma}$ be the normalization of
$\mul{\Sigma}$, so that $\Def_\Gamma(\mul{\Sigma})$ is isomorphic to
the space of deformations of $\ti{\Sigma}$ equipped with the
additional markings $w_1^\pm, \ldots, w_m^\pm$ obtained by lifting the
nodes.  The general theory of deformations, see for example
\cite{dou:ana} in the analytic setting, shows that any marked nodal
curve $\mul{\Sigma}$ admits a versal deformation with smooth parameter
space $S$.  $\mul{\Sigma}$ admits a universal deformation
$\mul{\Sigma}_S \to S$ if and only if $\mul{\Sigma}$ is stable.
Furthermore, the space $\Def(\mul{\Sigma})$ of the space of
infinitesimal deformations admits a canonical isomorphism with
$H^{0,1} ( \mul{\Sigma}, T \mul{\Sigma}[ - z_1 - \ldots - z_n ]) ,$
where $T\mul{\Sigma}[-z_1-\dots, -z_n]$ is the sheaf of vector fields
vanishing at $z_1, \dots, z_n$.

The relationship between the various deformation spaces (in the case
with markings, fixed type, etc.) is given as follows.  The space of
{\em infinitesimal automorphisms} $\aut(\mul{\Sigma},\mul{z})$ of
$(\mul{\Sigma},\mul{z})$ is the space
$\Vect(\mul{\Sigma},\mul{z}) = H^0( \mul{\Sigma}, T \mul{\Sigma}[ - z_1 -
  \ldots - z_n ]) $
of holomorphic vector fields vanishing at the marked points.  The
short exact sequence of sheaves
$$ 0 \to 
\oplus_{i=1}^n T_{z_i} \mul{\Sigma} \to  T \mul{\Sigma} 
\to T \mul{\Sigma} [ - z_1 - \ldots - z_n ] \to 0$$
gives a long exact sequence in cohomology \cite[p. 94]{hm:mc}
$$ 0 \to \Vect(\mul{\Sigma},\mul{z}) \to  \Vect(\mul{\Sigma}) \to \bigoplus_{i=1}^n T_{z_i}
\mul{\Sigma} \to \Def(\mul{\Sigma},\mul{z}) \to \Def(\mul{\Sigma}) \to
0. $$
From now on, we omit the markings from the notation, and study
deformations of a nodal marked curve $\mul{\Sigma} = (\mul{\Sigma},
\mul{z})$.  By $T_{w_i^\pm} \mul{\Sigma}$, we mean the tangent space in
the component of $\mul{\Sigma}$ containing $w_i^\pm$.  A {\em gluing
  parameter} for the $j$-th node is an element $ \delta_i \in
T_{w_i^+} \mul{\Sigma} \otimes T_{w_i^-} \mul{\Sigma} .$ The canonical
conormal sequence \cite[p. 100]{hm:mc} gives rise to an exact sequence
\begin{equation} \label{type}
 0 \to \Def_\Gamma(\mul{\Sigma}) \to \Def(\mul{\Sigma}) \to
\bigoplus_{i=1}^m T_{w_i^+} \mul{\Sigma} \otimes T_{w_i^-} \mul{\Sigma}
\to 0 .\end{equation}
After trivialization of the tangent spaces the gluing parameters are
identified with complex numbers.

Universal deformations of a smooth marked curve can be constructed for
example using Teichm\"uller theory \cite{earle:fibre} or by Hilbert
scheme methods \cite[p. 102]{hm:mc}.  Later we will need an explicit
gluing construction of a universal deformation of a stable marked
curve.  This construction seems to be well-known, but the only proof
we could find in the literature is Siebert \cite{sie:gen}.  The idea
is to remove small neighborhoods of the nodes, and glue the remaining
components together.  A {\em local coordinate near a smooth point} $z
\in \mul{\Sigma}$ is a neighborhood $U$ of $z$ and a holomorphic
isomorphism $\kappa$ of $U$ with a neighborhood of $0$ in the tangent
line $T_{z} \mul{\Sigma}$, whose differential $T_0 U\to T_z \Sigma$ is
the identity.

\begin{remark} \label{convex}
The space of local coordinates near $z$ is convex, since if
$\kappa_0,\kappa_1$ are local coordinates then any combination $t
\kappa_0 + (1-t)\kappa_1 $ is still holomorphic and has the same
differential at $z$, and so by the inverse function theorem is a
holomorphic isomorphism in a neighborhood of $z$.
\end{remark} 
\noindent
Any gluing parameter $\delta_i$ induces an identification
$$ T_{w_i^+} \mul{\Sigma}- \{ 0 \} \to 
  T_{w_i^-} \mul{\Sigma} - \{ 0 \}, \ \ \lambda_i^+ \mapsto 
\delta_i/ \lambda_i^-  .$$
Given local coordinates for the nodes of $\mul{\Sigma}$ and a set of
gluing parameters $\mul{\delta} = (\delta_1,\ldots,\delta_m)$, define a
(possibly nodal) curve $\mul{\Sigma}^{\mul{\delta}}$ by gluing together
small disks around the node $w_i$ by $z \mapsto \delta_i/z$, for every
gluing parameter $\delta_i$ that is non-zero, where $z$ is the local
coordinate given by $\kappa_i$. That is,
\begin{equation} \label{glue1}
 \mul{\Sigma}^{\mul{\delta}} = \bigcup_{i=1}^k \Sigma_i - \{ w_1^\pm,
 \ldots, w_m^\pm \} / ( z \sim (\kappa_i^+)^{-1} (\delta_i / \kappa_i^-(z) )  , i =
 1,\ldots, m )
\end{equation}
for pairs of points in the two components such that both coordinates
are defined.  In particular, the choice of local coordinates
near the nodes defines a splitting of the sequence \eqref{type}. 

The gluing construction works in families as follows.  Let
$I^{i,\pm}_\Gamma \to S_\Gamma $ resp. $\mul{I}_\Gamma \to S_\Gamma$
denote the vector bundle whose fiber at $s \in S_\Gamma$ is the
tangent line at the $j$-node resp.  tensor product of tangent lines at
the nodes,
\begin{equation} \label{gluingpar}
 I^{i,\pm}_{\Gamma,s} = T_{w_{i,s}^\pm} \mul{\Sigma}_s, \quad \mul{I}_{\Gamma,s} = \bigoplus_{j=1}^m T_{w_{i,s}^-} \mul{\Sigma}_s \otimes
T_{w_{i,s}^+} \mul{\Sigma}_s .\end{equation} 
Let $\mul{\Sigma}_{S_\Gamma} \to S_\Gamma$ be a family of nodal curves
of the same combinatorial type $\Gamma$, with nodal points
$(w_{S_\Gamma,j}^\pm)_{i = 1}^m$.  A {\em holomorphic system of local
  coordinates for the $i$-th node} is a holomorphic map $\kappa_i$
from a neighborhood $U_{i,\pm}$ of the zero section in $ I^{i,\pm}_S$
to $\mul{\Sigma}_S$ which is an isomorphism onto its image and induces
the identity at any point in the zero section.  Given a holomorphic
system of coordinates for each node $\mul{\kappa} =
(\kappa_1^+,\kappa_1^-, \ldots, \kappa_m^+,\kappa_m^-)$ the gluing
construction \eqref{glue1} produces a family $ \mul{\Sigma}_S \to S $
over an open neighborhood $S$ of the zero section in the bundle $I \to
S_\Gamma$ of gluing parameters.

\begin{theorem} \label{aut}  
 \cite[Proposition 2.4]{sie:gen} If $\mul{\Sigma}_{\Gamma,S}$ is a
 family giving a universal deformation of fixed type, then
 $\mul{\Sigma}_S $ is a universal deformation of any of its fibers, and
 in particular is independent up to isomorphism of deformations of the
 choice of local coordinates $\kappa$.
\end{theorem} 

\noindent The following properties of universal deformations of stable
curves will be used later:

\begin{theorem} \label{props}  \cite[Lemma 2.7]{sie:gen}
For any universal deformation $\mul{\Sigma}_S$, the action of
automorphisms $\Aut(\mul{\Sigma})$ of $\mul{\Sigma}$ extends to an
action of $\Aut(\mul{\Sigma})$ on $\mul{\Sigma}_S$, possibly after
shrinking $S$.  For any universal deformation, there exists a
neighborhood of the central fiber such that any two fibers
$\mul{\Sigma}_S$ contained in the neighborhood are isomorphic, if
and only if they are related by an automorphism of $\mul{\Sigma}$.
\end{theorem}


If $\mul{\Sigma}$ is not stable, then the above construction produces a
{\em minimal versal deformation} of $\mul{\Sigma}$.  That is,
$\mul{\Sigma}_S \to S$ is versal, and any other versal deformation
given by a family $\mul{\Sigma}'_{S'} \to S'$ is obtained by pull-back
by a map $S' \to S$.

Algebraic families of connected stable nodal curves with genus $g$ and
$n$ markings form the objects of a smooth {\em Deligne-Mumford stack}
$\ol{M}_{g,n}$ \cite{dm:irr} which admits a coarse moduli space with
the structure of a normal projective variety.  The maps
$\Def(\mul{\Sigma}) \to \ol{M}_{g,n}, \ s \mapsto [\mul{\Sigma}_s] $
(restricted to a neighborhood of $0$) provide $\ol{M}_{g,n}$ with an
atlas of holomorphic orbifold charts.

\subsection{Stratified-smooth families of stable curves} 
\label{parsmooth}

We extend the definition of families and deformations to smooth and
stratified-smooth settings.  Given a family $\mul{\Sigma}_S \to S$ of
compact complex nodal curves, let
$$S = \bigcup S_\Gamma,\quad S_\Gamma = \{ s \in S,
\ \ \Gamma(\mul{\Sigma}_s) = \Gamma \} $$ 
denote the stratification by combinatorial type of the fiber.  It
follows from the gluing construction of the previous section that if
$\mul{\Sigma}_S \to S$ is a family giving a universal deformation, then
each $S_\Gamma$ is a smooth manifold, and the restriction
$\mul{\Sigma}_{\Gamma,S_\Gamma}$ of $\mul{\Sigma}_{S_\Gamma}$ to
$S_\Gamma$ gives a universal deformation of fixed type $\Gamma$.  By a
{\em smooth family} of curves of fixed type $\Gamma$ we mean a fiber
bundle $\mul{\Sigma}_{\Gamma,S_\Gamma} \to S_\Gamma$ with fibers of
type $\Gamma$ and smoothly varying complex structure.  In the nodal
case, it is obtained from a smooth family of smooth holomorphic
curves, identified using a collection of pairs of smooth sections
(nodes).

\begin{lemma}  \label{smoothhol} 
 Holomorphic universal deformations of fixed type are also universal
 in the category of smooth deformations of $\mul{\Sigma}$.  That is,
 let $\mul{\Sigma}_S \to S , \varphi$ be a universal holomorphic
 deformation of fixed type of a nodal curve $\mul{\Sigma}$.  Any smooth
 deformation $\mul{\Sigma}_{S'}' \to S', \varphi'$ of nodal curves of
 fixed type is obtained by pull-back $\mul{\Sigma}_S \to S$ by a
 smooth map $S' \to S$.
\end{lemma} 

\begin{proof}   By the construction of local
 slices for the action of diffeomorphisms in \cite{earle:fibre},
 \cite[Chapter 9]{rob:const}.
\end{proof}  

Similarly we can define {\em continuous} families of holomorphic
curves, which correspond to continuous maps $S' \to S$ to the
parameter space $S$ for a universal holomorphic deformation.  The
following spells out the definition without reference to the universal
holomorphic deformation.

\begin{definition}  Let $\Gamma_0,\Gamma_1$ be graphs. 
A {\em simple contraction} $\tau$ is a pair of maps $\on{Vert}(\tau):
\on{Vert}(\Gamma_0) \to \on{Vert}(\Gamma_1)$ and a bijection
$\Edge(\tau): \Edge(\Gamma_0) \to \Edge(\Gamma_1) \cup \{ \emptyset
\}$ such that $\Gamma_1$ is obtained from $\Gamma_0$ identifying the
head and tail of the {\em contracting edge} $e$ such that
$\Edge(\tau)(e) = \emptyset$. A {\em contraction} is a sequence of
simple contractions.
\end{definition}

\begin{definition} \label{familycurves} A {\em continuous family} of nodal holomorphic curves
 consists of topological spaces $\mul{\Sigma}_S$, a surjection
 $\mul{\Sigma}_S \to S$, and a collection of (possibly nodal)
 holomorphic structures $j_{{\mul{\Sigma}}_s}$ on the fibers
 $\mul{\Sigma}_s, s \in S$, which vary continuously in $s$ in the
 following sense: for every $s_0 \in S$ there exists for $s$ in a
 neighborhood of $s_0$ of some combinatorial type $\Gamma$,
\begin{enumerate}
\item contractions $\tau_s : \Gamma(\mul{\Sigma}_{s_0}) \to
  \Gamma $, constant in $s \in S_\Gamma$;
\item for every node $ \{ w_i^\pm \}$ collapsed under $\tau_s$, a pair
  of local coordinates $\kappa^\pm_i:  W^\pm_i \to \C$ 
\item for every component $\Sigma_{s_0,i}$ of $\mul{\Sigma}_{s_0}$, a
  real number $\eps_s > 0$ and maps
$$\phi_{i,s}: \Sigma_{s_0,i} - \cup_{w_k^\pm \in \Sigma_{s_0,i},
    \tau_s(w_k^\pm) = \emptyset} (\kappa_k^\pm)^{-1}(B_{\eps_s}) \to
  \mul{\Sigma}_{s,\tau_s(i)}$$
\end{enumerate} 
such that
\begin{enumerate}
\item for any $s$, the images of the maps $\phi_{i,s}$
cover $\mul{\Sigma}_{s}$;
\item for any nodal point $w_i^\pm$ of $\mul{\Sigma}_s$ joining
  components $\Sigma_{s,i^\pm(k)}$, there exists a constant $\lambda_s
  \in \C^*$ such that $(\kappa_k^+ \circ \phi_{s,i^+(k)}^{-1} \circ
  \phi_{s,i^-(k)} \circ (\kappa_k^-)^{-1})(z) = \lambda_s z $, if the
  former is defined, and $\lambda_s \to 0$ as $s \to s_0$.
\item for any $z \in \Sigma_{s_0,i}$ in the complement of the
  $W_{k,s}^\pm$, $\lim_{s \to s_0}(\phi_{i,s}(z)) = z$;
\item $\phi_{i,s}^* j_{\Sigma_{s,\tau_s(i)}}$ converges to $
  j_{\Sigma_{s_0,i}}$ uniformly in all derivatives on compact sets;
\item if $z_i$ is contained in $\Sigma_{s_0,k}$, then $z_i = \lim_{s
  \to s_0} \phi_{s,k}^{-1}(z_{i,s})$.
\end{enumerate}
A {\em stratified-smooth} family of curves is a continuous family
$\mul{\Sigma}_S \to S$ over a stratified base $S =
\bigcup_{\Gamma} S_\Gamma$ such that the restriction
$\mul{\Sigma}_{S_\Gamma}$ of $\mul{\Sigma}_S$ to $S_\Gamma$
is a smooth family of fixed type $\Gamma$. 
A {\em stratified-smooth deformation} of a nodal curve
  $\mul{\Sigma}$ is a germ of a stratified-smooth family of nodal
  curves $\mul{\Sigma}_S$ equipped with an isomorphism of the
  central fiber $\mul{\Sigma}_0$ with $\mul{\Sigma}$.
A {\em universal stratified-smooth deformation} of $\mul{\Sigma}$
  is a deformation with the property that any other stratified-smooth
  deformation $\mul{\Sigma}_{S'}' \to S'$ is obtained by pull-back
  by maps $\psi: S' \to S$, $\phi: \mul{\Sigma} \times_S S' \to
  \mul{\Sigma}'$, and any two isomorphisms $\phi,\phi'$ inducing the
  identity on $\mul{\Sigma}$ are equal.
\end{definition}   

Any universal holomorphic deformation is also a universal
stratified-smooth deformation, essentially by Lemma \ref{smoothhol}.
In the stratified-smooth setting, the analog of Theorem \ref{props}
fails and we need an additional definition:

\begin{definition} \label{excellent} A universal stratified-smooth deformation 
$(\pi: \mul{\Sigma}_S \to S,\phi)$ is {\em strongly universal} if $\pi$
  is a universal deformation of any of its fibers, and two fibers of
  $\pi$ are isomorphic, if and only if they are related by the action
  of $\Aut(\mul{\Sigma})$.
\end{definition} 

The construction of universal deformations extends to the smooth
setting as follows. Let $\mul{\Sigma}_{S_\Gamma} \to S_\Gamma$ be a
smooth family of curves of fixed type $\Gamma$.  A {\em smooth system
  of local coordinates for the $i$-th node of
  $\mul{\Sigma}_{S_\Gamma}$} is a smooth map $\kappa_i$ from a
neighborhood $U_{i,\pm}$ of the zero section in $ I^{i,\pm}$ to
$\mul{\Sigma}_{S_\Gamma}$ which is an isomorphism onto its image and
induces the identity at zero.  Given a universal deformation
$(\mul{\Sigma}_{S_\Gamma} \to S_\Gamma,\varphi)$ of fixed type $\Gamma$
and a smooth system of local coordinates, applying the gluing
construction \eqref{glue1} gives a smooth family $ \mul{\Sigma}_S \to
S$ over an open neighborhood $S$ of $0$ in $\Def(\mul{\Sigma})$.  We
may identify $S$ with $\Def(\mul{\Sigma})$, for simplicity of notation.

\begin{theorem}  \label{log}
Let $\mul{\Sigma}$ be a stable curve.  The family $\mul{\Sigma}_S \to  S
\subset  \Def(\mul{\Sigma})$ constructed by gluing from a family
$\mul{\Sigma}_{\Gamma,S} \to S \subset \Def_\Gamma(\mul{\Sigma})$ of
fixed type, using any smooth family of local coordinates $\mul{\kappa}$
near the nodes, gives a strongly universal stratified-smooth
deformation of $\mul{\Sigma}$.
\end{theorem} 

\begin{proof}
Let $\mul{\Sigma}_{S^{\mul{\kappa}}}^{\mul{\kappa}} \to S^{\mul{\kappa}}$
be a family constructed via gluing using a smooth family of local
coordinates $\mul{\kappa}$ as in \eqref{glue1}, and $\mul{\Sigma}_S \to
S$ a universal deformation using a holomorphic family of local
coordinates by the same construction \eqref{glue1}.  By universality,
there exists a map $\psi: S^{\mul{\kappa}} \to S$ so that
$\mul{\Sigma}_{\psi(s)} \cong \mul{\Sigma}_s^{\mul{\kappa}}$.  It
suffices to show that $\psi$ is a diffeomorphism on each stratum.
Consider the canonical map from $T_{\mul{\delta}} S^{\mul{\kappa}}$ to
$\Def(\mul{\Sigma}^{\mul{\delta}})$, which maps an infinitesimal change
in the parameter space $S^{\mul{\kappa}}$ to the corresponding
infinitesimal deformation of $\mul{\Sigma}^{\mul{\delta}}$, which we
identify with an element of $\Omega^0(\mul{\Sigma}, \End(T
\mul{\Sigma}^{\mul{\delta}}))$.  Let $U \subset
\mul{\Sigma}^{\mul{\delta}}$ denote the gluing region, that is, the
image of the union of domains of the local coordinates.  The
deformations generated by the gluing parameters are supported in the
gluing region $U$.  On the other hand, linearly independent
deformations of fixed type $\Def_\Gamma(\mul{\Sigma})$ generate
deformations of the glued curve that are linearly independent
on${\mul{\delta}} - U$, for sufficiently small $U$.  (The generated
deformations will not vanish on $U$, because of the varying local
coordinates.) Thus the map $\Def_\Gamma(\mul{\Sigma}) \to
\Omega^0(\mul{\Sigma} - U)$ is injective; it follows that
$TS^{\mul{\kappa}} \to \Def(\mul{\Sigma}^{\mul{\delta}})$ is injective,
hence an isomorphism by a dimension count.  This shows that the map
$S^{\mul{\kappa}} \to S$ is a covering.  Let $\mul{\kappa}_t$ be a
family of local coordinates interpolating between $\mul{\kappa}$ and a
holomorphic family.  The corresponding family $\psi_t$ interpolates
between the identity and $\psi$.  Since each $\psi_t$ is a covering
and $\psi_0$ is the identity, each $\psi_t$ is a diffeomorphism.
\end{proof} 

The strongly universal deformations above defined using smooth
families of local coordinates provide smooth orbifold charts on
$\ol{M}_{g,n}$.  Since the space of local coordinates is convex, one
can construct the local coordinates for each stratum compatibly.
Namely, let $\Gamma'$ be a combinatorial type degenerating to
$\Gamma$.  Local coordinates for the nodes of $M_{g,n,\Gamma}$
induce local coordinates for $M_{g,n,\Gamma'}$, in a neighborhood
of $M_{g,n,\Gamma}$, via the gluing construction \eqref{glue1}.
\begin{definition}  \label{compatcoord}
A {\em compatible system of local coordinates} for $\ol{M}_{g,n}$ is a
system of local coordinates for the nodes of each stratum
$M_{g,n,\Gamma}$, so that the local coordinates on any stratum
$M_{g,n,\Gamma'}$ are induced from those on
$M_{g,n,\Gamma}$, in a neighborhood of $M_{g,n,\Gamma}$.
\end{definition} 
\noindent Compatible systems of local coordinates can be constructed by
induction on the dimension of $M_{g,n,\Gamma}$, using convexity
on the space of local coordinates in Remark \ref{convex}.

One can modify the gluing construction above by choosing a different
smooth structure on the space of gluing parameters.  In the language
of Hofer, Wysocki and Zehnder \cite[Appendix]{ho:sc},
\begin{definition} \label{glueprof} A
{\em gluing profile} is a diffeomorphism $\varphi: (0,1] \to
  [0,\infty)$.  The diffeomorphism given by $\varphi(\delta) = -1 +
    1/\delta$ will be called the {\em standard gluing profile};
    $\varphi(\delta) = e^{1/\delta} - e $ will be called the {\em
      exponential gluing profile}, and $\varphi(\delta) = -\ln(\delta)
    $ the {\em logarithmic gluing profile}.
\end{definition} 
The set of gluing profiles naturally forms a partially ordered set:
Write $\varphi_1 \geq \varphi_0$ and say $\varphi_1$ is {\em softer
  than} $\varphi_0$ if $\varphi_1^{-1} \varphi_0$ extends to a
diffeomorphism of $[0,1]$. Write $\varphi_1 > \varphi_0$ and say that
$\varphi_1$ is {\em strictly softer} than $\varphi_0$ if the
derivatives of $\varphi_1^{-1} \varphi_0: [0,1] \to [0,1]$ vanish at
$0$.  The exponential gluing profile, standard gluing profile, and
logarithmic gluing profile form a decreasingly soft sequence in this
partial order.  

Fix a gluing profile $\varphi$, and consider once again the gluing
construction.

\begin{definition}  Given a nodal curve $\mul{\Sigma}$ with local coordinates
$\mul{\kappa}$ near the nodes, and a collection of gluing parameters
  $\delta = (\delta_1,\ldots,\delta_m)$, the {\em glued curve}
  $\mul{\Sigma}(\mul{\delta},\varphi)$ is defined by gluing together
  small disks:
\begin{equation} \label{gluecurve}
\mul{\Sigma}^{\mul{\delta},\varphi,\kappa} := \left( \bigcup_{i=1}^m
\Sigma_i - \{ w_1^\pm, \ldots, w_m^\pm \} \right) /
\sim \end{equation}
\noindent where the equivalence relation $\sim$ is given by 
$$ z \sim (\kappa_i^+)^{-1}( \exp( - \varphi(|\delta_i|) - i
\arg(\delta_i))/\kappa_i^-(z) , \quad, z \in U_i^-, \quad i =
1,\ldots, m .$$
More generally, given a family $\mul{\Sigma}_{S_\Gamma} \to S_\Gamma$
of curves of constant combinatorial type $\Gamma$ and a system of
local coordinates near the nodes $ \mul{\kappa}$, the construction
\eqref{gluecurve} produces a family of curves $
\mul{\Sigma}_{S^{\mul{\kappa},\varphi}}^{\varphi,\mul{\kappa}} \to
S^{\mul{\kappa},\varphi} $ where $S^{\mul{\kappa},\varphi}$ is the
product of $S$ with the space of gluing parameters.
\end{definition}
Let $\mul{\Sigma}$ be a compact, complex nodal curve.
For any gluing profile $\varphi$ and any collection $\mul{\kappa}$ of
local coordinates near the nodes, the family $
\mul{\Sigma}_{S^{\mul{\kappa},\varphi}}^{\mul{\kappa},\varphi} \to
S^{\mul{\kappa},\varphi} $ is a stratified-smooth strongly universal
deformation, since it is so for the standard gluing profile. 
Let $\mul{\kappa} = (\mul{\kappa}_\Gamma)$ be a compatible system of
local coordinates near the nodes, for each combinatorial type
$\Gamma$.  Each stratified-smooth universal deformation above defines
a {\em classifying map}
\begin{equation} \label{gluemaps}
S^{\mul{\kappa},\varphi} /\Aut(\mul{\Sigma}) \to \ol{M}_{g,n}, \quad s
\mapsto [\mul{\Sigma}_s] \end{equation}
which is a homeomorphism onto its image, possibly after shrinking the
parameter space $ S^{\mul{\kappa},\varphi} $.  (To obtain a precise
meaning for ``classifying map'' it is necessary to pass to the
stacks-theoretic viewpoint, which we do not discuss here.)  The maps
\eqref{gluemaps} provide $\ol{M}_{g,n}$ with a compatible set of
stratified-smooth orbifold charts, since the transition maps are the
identity on the space of gluing parameters by construction, and smooth
on each stratum.  We denote by $\ol{M}_{g,n}^{\mul{\kappa},\varphi}$
the smooth structure on $\ol{M}_{g,n}$ defined by the system of local
coordinates $\mul{\kappa}$ near the nodes and the gluing profile
$\varphi$; the use of this smooth structure seems to have been
suggested by Hofer.  It seems that these smooth structures might
depend on the choice of $\mul{\kappa}$, except in the case of the
logarithmic gluing profile, in which case one has a canonical smooth
structure.

The forgetful maps with respect to these non-standard smooth
structures have regularity properties that are worse than those with
respect to the standard smooth structure.  For $2g + n > 3$ we have
forgetful morphisms $f_i: \ol{M}_{g,n} \to \ol{M}_{g,n-1}$ by
forgetting the $i$-th marking and collapsing unstable components.
There are two possibilities: a genus zero component with one marking
and two nodes is replaced by a point; a genus zero component with two
markings and one node is replaced by a single marking.  For any gluing
profile, the maps $f_i$ are smooth away from the locus where
collapsing occurs.  We say a local coordinate on a genus zero curve is
{\em standard} if it extends to an isomorphism with the projective
line.  The forgetful morphism $f_i$ is smooth near the locus of one
node, two marking components if the local coordinates are standard and
$\delta \mapsto \exp(\varphi(\delta))^{-1}$ is smooth, that is,
$\varphi$ is at least as hard as the logarithmic gluing profile.  The
forgetful morphism $f_i$ is smooth near the locus of curves containing
components with two nodes and one marking if the map
$\delta_1,\delta_{2} \mapsto \varphi^{-1}( \varphi(\delta_1) +
\varphi(\delta_{2})) $ is smooth.  For example, in the logarithmic
gluing profile we have $(\delta_1, \delta_2) \mapsto \delta_1
\delta_2$, which is smooth, while for the standard gluing profile
collapsing a component gives the map $ (\delta_1,\delta_2) \mapsto
\delta_1 \delta_2 / (\delta_1 + \delta_2) $ in the local gluing
parameters, which is not smooth.

\section{Deformations of holomorphic maps from curves} 
\label{maps} 

This section reviews the construction of a stratified-smooth universal
deformations for stable (pseudo)holomorphic maps.  The proof relies on
a gluing theorem, of the sort given by Ruan-Tian \cite{rt:hi}; our
approach follows that of McDuff-Salamon \cite{ms:jh} who treat the
genus zero case.  A different set-up for gluing is described in
Fukaya-Oh-Ohta-Ono \cite{fooo}, and explained in more detail in
Abouzaid \cite{ab:ex}.  The gluing construction gives rise to charts
for the moduli space of regular stable maps.

\subsection{Stable maps} 

Let $(X,\omega)$ be a compact symplectic manifold and $\J(X)$ the
space of compatible almost complex structures on $X$.  Let $J \in
\J(X)$.

\begin{definition} 
A {\em marked nodal $J$-holomorphic map} to $X$ consists of a nodal
curve $\mul{\Sigma}$, a collection $\mul{z} = (z_1,\ldots, z_n)$ of
distinct, smooth points on $\mul{\Sigma}$, and a $J$-holomorphic map
$\mul{u}: \mul{\Sigma} \to X$.  
An {\em isomorphism} of marked nodal maps from $(\mul{\Sigma}_0,
\mul{z}_0, \mul{u}_0)$ to $(\mul{\Sigma}_1, \mul{z}_1, \mul{u}_1)$ is an
isomorphism of nodal curves $\mul{\psi}: \mul{\Sigma}_0 \to
\mul{\Sigma}_1$ such that $\mul{\psi}(z_{0,i}) = z_{1,i}$ for $i=1,
\ldots, n$ and $\mul{u}_1 \circ \mul{\psi} = \mul{u}_0$.  A marked nodal
map $(\mul{\Sigma},\mul{u},\mul{z})$ is {\em stable} if it has finite
automorphism group or equivalently each component $\Sigma_i$ of genus
zero resp. one for which $u_i$ is constant has at least three
resp. one special (nodal or marked) point.  The {\em homology class}
of stable map $\mul{u}: \mul{\Sigma} \to X$ is $\mul{u}_* [\mul{\Sigma}]
\in H_2(X,\Z)$.
\end{definition}  

A \emph{continuous family} of $J$-holomorphic maps over a topological
space $S$ is a continuous family of nodal curves $\mul{\Sigma}_S\to S$
(see Definition \ref{familycurves}) and a continuous map $\mul{u}:
\mul{\Sigma}_S\to X$ which is fiberwise holomorphic.  That is, for
each $s_0 \in S$ and each nearby combinatorial type $\Gamma$ we have
\begin{enumerate}
\item a sequence of contractions
$\tau_s : \Gamma(\mul{\Sigma}_{s_0}) \to \Gamma $, constant in $s \in
  S_\Gamma$;
\item for every node $ \{ w_i^\pm \}$ collapsed under $\tau_s$, a pair
  of local coordinates $\kappa^\pm_i: W^\pm_i \to \C$
\item for every component $\Sigma_{s_0,i}$ of $\mul{\Sigma}_{s_0}$, a
  real number $\eps_s > 0$ converging to $0$ as $s \to s_0$ and maps
$$\phi_{i,s}: \Sigma_{s_0,i} - \cup_{w_k^\pm \in \Sigma_{s_0,i},
    \tau_s(w_k^\pm) = \emptyset} (\kappa_k^\pm)^{-1}(B_{\eps_s}) \to
  \mul{\Sigma}_{s,\tau_s(i)}$$
\end{enumerate} 
such that
\begin{enumerate}
\item for any $s$, the images of the maps $\phi_{i,s}$
cover $\mul{\Sigma}_{s}$;
\item for any nodal point $w_i^\pm$ of $\mul{\Sigma}_s$ joining
  components $\Sigma_{s,i^\pm(k)}$, there exists a constant $\lambda_s
  \in \C^*$ such that $(\kappa_k^+ \circ \phi_{s,i^+(k)}^{-1} \circ
  \phi_{s,i^-(k)} \circ (\kappa_k^-)^{-1})(z) = \lambda_s z $ where
  defined, and $\lambda_s \to 0$ as $s \to s_0$.
\item for any $z \in \Sigma_{s_0,i}$ in the complement of the
  $W_{k,s}^\pm$, $\lim_{s \to s_0}(\phi_{i,s}(z)) = z$;
\item $\phi_{i,s}^* j_{\Sigma_{s,\tau_s(i)}}$ converges to $
  j_{\Sigma_{s_0,i}}$ uniformly in all derivatives on compact sets;
\item if $z_i$ is contained in $\Sigma_{s_0,k}$, then $z_i = \lim_{s
  \to s_0} \phi_{s,k}^{-1}(z_{i,s})$.
\item $\phi_{i,s}^* u_s$ converges to $u_{s_0}$ uniformly in all
  derivatives on compact sets.
\end{enumerate}
\begin{remark} It follows from the assumption that $u_S: \Sigma_S \to X$
is continuous that the homology class $u_{s,*} [\Sigma_s]$ is locally
constant in $s \in S$. Indeed continuity implies that for $s$
sufficiently close to $s_0$, $u_S$ is homotopic to a map of the form
$v_S \circ \gamma_S$ where $\gamma_s: \Sigma_S \to \Sigma_{s_0}$ is a
map to the central fiber $\Sigma_{s_0}$ which collapses the gluing
regions to the node.  Since each $\gamma_s = \gamma_S | \Sigma_s$ maps
$[\Sigma_s]$ to $[\Sigma_{s_0}]$, the claim follows.
\end{remark} 

In particular, taking $S$ to be the topological space given as the
closure of the set $S^*$ of rational numbers of the form $1/i, i \in
\Z_{> 0}$, we say that a sequence of holomorphic maps $u_i: \Sigma_i
\to X$ {\em Gromov converges} if it extends to a continuous family
over $S$.  To state the Gromov compactness theorem, recall that the
{\em energy} of a map $u: \Sigma \to X$ is
$$ E(u)=\frac{1}{2}\int |du|^2 .$$

\begin{theorem} [Gromov compactness] 
Let $X,\omega,J$ be as above.  Any sequence $u_i: \Sigma_i \to X$ of
stable holomorphic maps with bounded energy has a Gromov convergent
subsequence.  Furthermore, the limit is unique.
\end{theorem} 

For references and discussion, see for example \cite[Theorem
  1.8]{io:rel}.  The definition of Gromov convergence passes naturally
to equivalence classes of stable maps.  A subset $C$ of
$\ol{M}_{g,n}(X,d)$ is {\em Gromov closed} if any sequence in $C$ has
a limit point in $C$, and {\em Gromov open} if its complement is
closed.  The Gromov open sets form a topology for which any convergent
sequence is Gromov convergent, by an argument using \cite[Lemma
  5.6.5]{ms:jh}.  Furthermore, any convergent sequence has a unique
limit.  Gromov compactness implies that for any $E> 0$, the union of
$\ol{M}_{g,n}(X,d)$ over $d \in H_2(X,\Z)$ with $(d, [\omega]) < E$ is
a compact, Hausdorff space.

\begin{definition} \label{ssfam}
 Let $X,\omega,J$ be as above.
A {\em stratified-smooth} family of nodal $J$-holomorphic maps over a
space $S$ is a pair $(\mul{\Sigma}_S,\mul{u}_S)$ of a stratified-smooth
family of nodal curves $\mul{\Sigma}_S \to S$ together with a
continuous map $\mul{u}_S: \mul{\Sigma}_S \to X$ such that the
restriction $\mul{u}_s$ of $\mul{u}$ to any fiber $\mul{\Sigma}_s$ is
holomorphic, and the restriction of $\mul{u}_S$ to any stratum
$\mul{\Sigma}_\Gamma$ is smooth.
A {\em stratified-smooth deformation} of a stable $J$-holomorphic map
$(\mul{\Sigma}, \mul{u})$ is a germ of a stratified-smooth family
$(\mul{\Sigma}_S, \mul{u}_S)$ together with an isomorphism
of nodal maps $\iota: \mul{\Sigma}_0 \to \mul{\Sigma}$ such that
$\iota^* \mul{u} = \mul{u}_0$.
A deformation $(\mul{\Sigma}_S, \mul{u}_S,\iota)$ of
  $(\mul{\Sigma},\mul{u})$ is {\em versal} if any other (germ of) family
  of marked, nodal curves $(\mul{\Sigma}',\mul{\Sigma}_0) \to (S',0)$ is
  induced from a map $\psi: S' \to S$ in the sense that there exists
  an isomorphism $\phi: \mul{\Sigma}' \to \mul{\Sigma} \times_S S'$ in a
  neighborhood of the central fiber $ \mul{\Sigma}_0$, and $\mul{u}'$ is
  obtained by composing projection on the first factor with $\mul{u}$.
  A versal deformation is {\em universal} if the map $\phi$ above is
  the unique map inducing the identity on $\mul{\Sigma}_0$.
\end{definition}  

\subsection{Smooth universal deformations of regular stable maps of  fixed
combinatorial type.} 

Let $\mul{u} :\mul{\Sigma} \to X$ be a stable map.  For $p >2$ define
a fiber bundle $\E \to \B$ by
$$ \B = \J(\mul{\Sigma}) \times \Map(\mul{\Sigma},X)_{1,p} , \quad
\E_{\zeta,\mul{u}} = \Omega^{0,1}(\mul{\Sigma}, \mul{u}^* TX)_{0,p}
,$$
where the latter is the space of $(0,1)$-forms with respect to the
pair $(\mul{j}(\zeta),J)$.  Consider the Cauchy-Riemann section,
$$ \ol{\partial}:\B \to \E, \ \ \ (j,\mul{u}) \mapsto \olp_j \mul{u},
\ \ \olp_j \mul{u} = \hh ( \d \mul{u} \circ \mul{j}(\zeta) - J_u \circ
\d \mul{u}) .$$
Let
$$ \ev: \B \to X^{2m} , \ \ \ \mul{u} \mapsto 
(\mul{u}(w_1^-), \mul{u}(w_1^+), \ldots, 
\mul{u}(w_m^-), \mul{u}(w_m^+) ) $$ 
denote the map evaluating at the nodal points.  The space of stable
maps of type $\Gamma$ is given as $(\olp,\ev)^{-1}(0)$.  To obtain a
Fredholm map, we quotient by diffeomorphisms of $\Sigma$, or
equivalently, restrict to a minimal versal deformation $\mul{\Sigma}_S
\to S$ of $\mul{\Sigma}$ of fixed type.  This means that for each
$\zeta \in \Def_\Gamma(\mul{\Sigma})$ near $0$ we have a complex
structure $\mul{j}(\zeta)$ on $\mul{\Sigma}$, which we may assume
agrees with $\mul{j} = \mul{j}(0)$ near the nodes.  Then the
Cauchy-Riemann section induces a map
$$ \Def_\Gamma(\Sigma) \times \Omega^0(\Sigma, u^* TX) \to \cE .$$
Linearizing the Cauchy-Riemann section, together with the differences
at the nodes, gives rise to a Fredholm operator
\begin{equation} 
 \label{explinear} 
 \ti{D}_{\mul{u}} : \Def_\Gamma(\mul{\Sigma}) \times
 \Omega^0(\mul{\Sigma}, \mul{u}^* TX) \to
 \Omega^{0,1}(\mul{\Sigma},\mul{u}^* TX) \oplus \bigoplus_{i=1}^m
 \mul{u}(w_i^\pm)^* TX $$ $$ \ti{D}_{\mul{u}} (\mul{\zeta},\mul{\xi}) :=
 \left( \pi^{0,1}_{\mul{\Sigma}} ( \nabla \mul{\xi} - \hh J(\mul{u}) \d
 \mul{u} Dj(\mul{\zeta}) - \hh J_{\mul{u}} (\nabla_{\mul{\xi}} J )_{\mul{u}}
 \partial \mul{u})), ( \mul{\xi}(w_i^+)- \mul{\xi}(w_i^-) )_{i=1}^m
 \right)
\end{equation}
given by the linearized Cauchy-Riemann operator on each component, and
the difference of the values of the section at the nodes
$w^\pm_1,\ldots, w^\pm_m$.  The map $\mul{u} =
(\mul{\Sigma},\mul{u},\mul{z})$ is {\em regular} if $\ti{D}_{\mul{u}}$
is surjective.  This is independent of the choice of representatives
$j(\zeta)$: any two such choices $\mul{j}'(\zeta),\mul{j}(\zeta)$ are
related by a diffeomorphism of $\mul{\Sigma}$.  The {\em space of
  infinitesimal deformations of $\mul{u}$ of fixed type} is
$$ \Def_\Gamma(\mul{u}) = \ker(\ti{D}_{\mul{u}})/\aut(\mul{\Sigma}) .$$
 The {\em space of infinitesimal deformations of $\mul{u}$} is
$$ \Def(\mul{u}) = \Def_\Gamma(\mul{u}) \oplus 
\bigoplus_{i=1}^m
 T_{w_j^+} \mul{\Sigma} \otimes T_{w_j^-} \mul{\Sigma} $$
where $\Gamma$ is the type of $\mul{u}$.

\begin{theorem} \label{fixed}  Let $X,\omega,J$ be as above.
A regular marked nodal $J$-holomorphic map $\mul{u} =
(\mul{\Sigma},\mul{u},\mul{z})$ admits a strongly universal deformation
$(\mul{\Sigma}_S,\mul{u}_S,\mul{z}_S)$ with parameter space $S \subset
\Def_\Gamma(\mul{u})$ of fixed type if and only if $\mul{u}$ is stable.
\end{theorem}    


\begin{proof}  Let $(\mul{\Sigma},\mul{u})$ be a stable map to $X$ and
$\mul{\Sigma}_S \to S \subset \Def_\Gamma(\mul{\Sigma})$ a minimal
  versal deformation of $\mul{\Sigma}$ of fixed type constructed in
  \eqref{glue1}.  We may write any map $C^0$-close to $\mul{u}$ as
  $\exp_{\mul{u}} (\mul{\xi})$ for some $\mul{\xi} \in
  \Omega^0(\mul{\Sigma}, \mul{u}^* TX)$.  Let $\Psi_{\mul{u}}(\mul{\xi}):
  \mul{u}^* TX \to \exp_{\mul{u}}(\mul{\xi})^* TX$ denote parallel
  transport along geodesics with respect to the Hermitian connection $
  \ti{\nabla} = \nabla - \hh J(\nabla J) ;$ here $\nabla$ is the
  Levi-Civita connection, see \cite[Chapter 2]{ms:jh}.  This defines an
  isomorphism
\begin{equation} \label{parallel} 
\Psi_{\mul{u}}(\mul{\xi})^{-1} : \Omega^{0,1}_{\mul{j}}( \mul{\Sigma},
\exp_{\mul{u}}(\mul{\xi})^* TX) \to
\Omega^{0,1}_{\mul{j}}(\mul{\Sigma}, \mul{u}^* TX)
.\end{equation} 
where subscript $j$ denotes the space of $0,1$-forms taken with
respect to the complex structure $j$ on $\mul{\Sigma}$.  There is an
isomorphism of $\Omega^{0,1}_{\mul{j}(\zeta)}(\mul{\Sigma}, \mul{u}^*
TX)$ with $\Omega^{0,1}_{\mul{j}}(\mul{\Sigma},\mul{u}^*TX)$ given by
composing the inclusion
$$\Omega^{0,1}_{\mul{j}(\zeta)}(\mul{\Sigma},\mul{u}^* TX) \to \Omega^1(\mul{\Sigma},
\mul{u}^* TX)_\C = \Omega^1(\mul{\Sigma}; \mul{u}^* TX) \otimes_\R \C$$
with the projection $ \Omega^1(\mul{\Sigma}; \mul{u}^* TX) \otimes_\R \C
\to\Omega^{0,1}_{\mul{j}}(\mul{\Sigma},\mul{u}^* TX) .$ We denote by
$$ \Psi_j(\zeta): \Omega^{0,1}_{\mul{j}}(\mul{\Sigma}, \mul{u}^* TX) \to
\Omega^{0,1}_{\mul{j}(\zeta)}(\mul{\Sigma},\mul{u}^* TX) $$
the resulting map; one can think of this as a connection over the
space of complex structures on $\mul{\Sigma}$ on the bundle whose fiber
is the space of $0,1$-forms with respect to $j(\zeta)$. By composing
$\Psi_{\mul{u}}(\mul{\xi})^{-1}$ and $\Psi_j(\zeta)^{-1}$ we
obtain an identification
\begin{equation} \label{eq:para}
\Psi_{j,\mul{u}}(\zeta,\xi)^{-1} :
\Omega^{0,1}_{\mul{j}(\zeta)}(\mul{\Sigma},\exp_{\mul{u}}(\mul{\xi})^* TX) \to
\Omega^{0,1}_{\mul{j}}(\mul{\Sigma},\mul{u}^*TX) .
\end{equation}
Define
$$ \cF_{\mul{u}}: \Def_\Gamma(\mul{\Sigma}) \times
\Omega^0(\mul{\Sigma},\mul{u}^* TX) \to
\Omega^{0,1}_{\mul{j}}(\mul{\Sigma},\mul{u}^* TX) $$
$$ (\zeta,\mul{\xi}) \mapsto
\Psi_{j,\mul{u}}(\zeta,\mul{\xi})^{-1}
(\olp_{\mul{j}(\zeta)}(\exp_{\mul{u}}(\mul{\xi}))) .$$
The operator $\ti{D}_{\mul{u}}$ is the linearization of $\cF_{\mul{u}}$.
The implicit function theorem implies that if $\mul{u}$ is regular then
the zero set of $\cF_{\mul{u}}$ is modelled locally on a neighborhood
of $0$ in $\ker(\ti{D}_{\mul{u}})$.  Furthermore, by elliptic
regularity the zero set consists entirely of smooth $J$-holomorphic
maps \cite[Section B.4]{ms:jh}.  Thus we obtain a smooth family of
stable maps in a neighborhood of $0$ in $\ker(\ti{D}_{\mul{u}})$.  The
action of $\Aut(\mul{u})$ on the space of stable maps with domain
$\mul{\Sigma}$ induces an inclusion of the Lie algebra $\aut(\mul{u})$
into $\ker(\ti{D}_{\mul{u}})$.  Restricting to $\Def_\Gamma(\mul{u})$,
identified with a complement of $\aut(\mul{\Sigma})$ (that is, a slice
for the $\Aut(\mul{u})$ action) gives a family $(\mul{\Sigma}_S,
\mul{u}_S) \to S \subset \Def_\Gamma(\mul{u})$ of fixed type.  The
family $(\mul{\Sigma}_S,\mul{u}_S)$, together with the canonical
identification $\iota$ of the central fiber with $\mul{\Sigma}$, is a
universal smooth deformation of fixed type.  Indeed, another smooth
family $(\mul{\Sigma}_{S'},\mul{u}_{S'})$ over a base $S'$ is in
particular a deformation of the underlying curve. After shrinking
$S'$, each fiber of $(\mul{\Sigma}_{s'}, \mul{u}_{s'}')$ corresponds to
a zero of $\cF_{\mul{u}}$, and so lies in the image of the map given by
the implicit function theorem.  The uniqueness part of the implicit
function theorem gives a smooth map $\psi: S' \to \Def_\Gamma(\mul{u})$
and an identification $\mul{\Sigma}_{S'} \to \psi^* \mul{\Sigma}_S$.
Any two such maps inducing the same map on the central fiber are close
in a neighborhood of the central fiber.  Since the automorphism groups
of the central fibre are discrete, any automorphism group is
discrete. Thus any two such automorphisms defined in a neighborhood of
the central fiber, and equal on the central fiber must be equal in a
neighborhood of the central fiber.  This shows that the identification
is unique, so that the deformation given by the gluing construction is
universal.

If $\mul{u}$ is not stable, then it has no universal deformation since
the identification with the central fiber is unique only up to a
continuous family of automorphisms.
\end{proof} 

Let $M_{g,n,\Gamma}^{\reg}(X,d)$ denote the moduli space of regular
stable maps of combinatorial type $\Gamma$.  A family $ \mul{u}_S$ over
$S \subset \Def_\Gamma(\mul{u}) $ induces a map
\begin{equation} \label{fixedcharts}
S \to M_{g,n,\Gamma}^{\reg}(X,d), \ s \mapsto [\mul{u}_s]
\end{equation}
where $[\mul{u}_s]$ denotes the isomorphism class of $\mul{u}_s:
\mul{\Sigma}_s \to X$.

\begin{theorem} For any $g,n,d$ and combinatorial type $\Gamma$
with $m$ nodes, $M_{g,n,\Gamma}^{\reg}(X,d)$ has the structure of a
smooth orbifold of dimension $(1-g)( \dim(X) - 6) + 2 (c_1 (TX),d) -
2m + 2n$, with tangent space at $[\mul{u}]$ isomorphic to
$\Def_\Gamma(\mul{u})$.
\end{theorem}

\begin{proof}   
By Theorem \ref{fixed}, the maps \eqref{fixedcharts} for families
giving universal deformations are homeomorphisms onto their image and
provide compatible charts.  The dimension formula follows from
Riemann-Roch: The index of $\ti{D}_{\mul{u}}$ may be deformed to a
complex linear operator by homotoping the zero-th order terms (which
define a compact operator) to zero.
\end{proof}

\subsection{Constructing stratified-smooth deformations of varying type} 
\label{gluing}

The main result of this section is Theorem \ref{premain}, which is
probably well-known, cf.  \cite{rrs:mod}, \cite{rt:hi}, but for which
we could not find an explicit reference.  The theorem itself will not
be used, but the estimates involved in the proof will be needed later
for the corresponding result for vortices.  The proof uses a gluing
construction for holomorphic maps, which produces from a smooth family
of holomorphic maps of fixed type, a stratified-smooth family of maps
of varying type.

\noindent {\em Step 1:  Approximate Solution} 

\begin{definition} Let $\mul{\Sigma}$ be a compact, complex nodal curve. 
A {\em gluing datum} for $\mul{\Sigma}$ consists of 
\begin{enumerate} 
\item a collection of gluing parameters $\mul{\delta} =
  (\delta_1,\ldots,\delta_m)$ in the bundle $\mul{I}$ of
  \eqref{gluingpar};
\item local coordinates $\kappa_j^\pm$ near the nodes $w_j^\pm$ for
  $j=1,\ldots,m$;
\item a parameter $\rho$ which describes the width of the annulus on
  which the gluing of maps is performed;
\item a gluing profile $\varphi$, see Definition \ref{glueprof};
\item a smooth cutoff function
\begin{equation} 
\label{firstcut}
 \alpha: \C \to [0,1 ], \ \ \ \alpha(z) = \left\{ 
\begin{array}{ll} 0 & |z| \leq 1 \\
1 & |z| \geq 2 \end{array} \right\} .\end{equation} 
\end{enumerate}
\end{definition} 
We first treat the case that $\varphi$ is the standard gluing profile.
Let a gluing datum be given, and $ \mul{\Sigma}^{\mul{\delta}}$ denote
the glued curve from \eqref{glue1}.  Let $\mul{u}: \mul{\Sigma} \to X$
be a holomorphic map.  Near each node $w_k$ let $i^\pm(k)$ denote the
components on either side of $w_k$.  In the neighborhoods $U_k^\pm$
(assuming they have been chosen sufficiently small) define maps
$$ \xi_k^\pm: U_k^\pm \to T_{x_k} X, \ \ \ u_{i^\pm(k)}(z) =
\exp_{x_k}(\xi_k^\pm(z)) $$
where $x_k = u(w_k)$ and $\exp_{x_k}: T_{x_k} X \to X$ denotes
geodesic exponentiation.  Given a holomorphic map $\mul{u}:\mul{\Sigma}
\to X$, and a gluing datum
$(\mul{\delta},\mul{\kappa},\rho,\varphi,\alpha)$ define the {\em
  pre-glued map} by interpolating between the maps on the various
components using the given cutoff function and local coordinates:
 $\mul{u}^{\mul{\delta}} = \mul{u}(z)$ for $z \notin \cup_k U_k^\pm$ and 
otherwise 
\begin{equation} \label{preglued}
\mul{u}^{\mul{\delta}} (z) = \exp_{x_k}( \alpha(\kappa_k^\pm( z) / \rho
|\delta_k|^{1/2} ) \xi_k^\pm(z) ) \quad z \in U_k^\pm .\end{equation}

\begin{remark} \label{nodalrem} 
The same formula but with domain $\mul{\Sigma}$ (not the glued curve)
defines an {\em intermediate map} $ \mul{u}_0^{\mul{\delta}}:
\mul{\Sigma} \to X $ which is constant near the nodes.  The right
inverse of $\ti{D}_{\mul{u}_0^{\mul{\delta}}}$ will be used in the
gluing construction. \end{remark}

\begin{figure}
\includegraphics{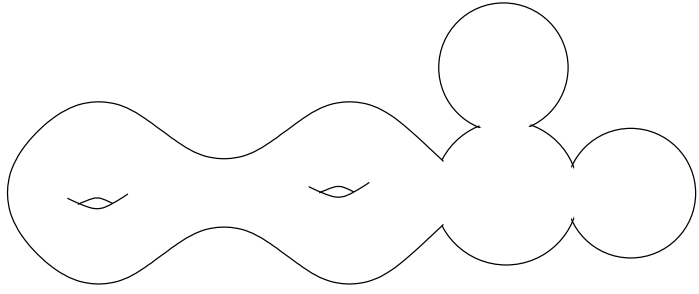}
\caption{Continuous metric on a glued curve}
\label{metric}
\end{figure}

First we estimate the failure of
$\mul{u}^{\mul{\delta}}$ to satisfy the Cauchy-Riemann equation.  Define
on $\mul{\Sigma}^{\mul{\delta}}$ the $C^0$-metric $g$ by the
identification
\begin{equation} \label{gluemetric} \mul{\Sigma}^{\mul{\delta}} = \mul{\Sigma} - \bigcup_{k,\pm}
\kappa_k^\pm(B_{|\delta_k|^{1/2}}(0)) / \left( \kappa^+_k(\partial B_{|\delta_k|^{1/2}}(0))\sim \kappa^-_k(\partial B_{|\delta_k|^{1/2}}(0)) \right)
\end{equation}
using a K\"ahler metric on $\Sigma$, see Figure \ref{metric}. 
\noindent The generalized Sobolev spaces $W^{l,p}$ with respect to
this metric are defined for $p \ge 1$ and integers $l \in \{0,1\}$,
see \cite{ad:so} or \cite{au:non}.  For any vector bundle $E$ we
denote by $\Omega(\mul{\Sigma}^{\mul{\delta}}, E)_{l,p,\mul{\delta}}$ the
space of $W^{l,p}$ forms with values in $E$.  If $p = \infty$ the norm
is independent of $\mul{\delta}$ and we drop it from the notation.  Let
$\Vert \cdot \Vert_{k,p,\mul{\delta}}$ denote the Sobolev
$W^{k,p}$-norm on $\Omega^0(\mul{u}^* TX)$ defined using the
$\mul{\delta}$-dependent metric \eqref{gluemetric}.

\begin{proposition} \label{errorprop}
Suppose that $\mul{u}: \mul{\Sigma} \to X$ is a stable map, and
$\mul{u}^{\mul{\delta}}: \mul{\Sigma}^{\mul{\delta}} \to X$ is the
pre-glued map defined in \eqref{preglued}, defined for $\mul{\delta}$
sufficiently small.  There is a constant $c$ and an $\eps > 0$ such
that if
$ \Vert \mul{\delta} \Vert < \eps,  
\rho > 1/\eps,$ and $ |\delta_k|^2 \rho < \eps, 
 \ \ \  k= 1,\ldots, m $
then 
$$\Vert \ol{\partial} \mul{u}^{\mul{\delta}} \Vert_{0,p,\mul{\delta}}^p
\leq c \sum_{k=1}^m (|\delta_k|^{1/2} \rho)^{2} .$$
\end{proposition}

\begin{proof}   
Compare with McDuff-Salamon \cite[Chapter 10]{ms:jh}. The error term
$\ol{\partial} \mul{u}(\mul{\delta})$ can be estimated by terms of two
types; those involving derivatives of the cutoff functions and those
involving derivatives of the map $\xi_k$.
%
The derivative of $\exp_{x_k}$ is approximately the identity near the
node.  The derivatives of $\alpha$ grow like $ 1/ \rho |
\delta_k|^{1/2}$, while the norm of $\xi_k^\pm$ is bounded by a
constant times $\rho |\delta_k |^{1/2}$ on the gluing region.  The
term involving the derivatives of $\alpha$ is bounded and supported
on a region of area less than $ \pi \rho^2 |\delta_k|$ for each
node. The derivatives of $\xi_k^\pm$ are also uniformly bounded, and
the area bound gives the required estimate.
\end{proof}  

Let $\mul{\Sigma}_S \to S$ with $S \subset \Def_\Gamma(\mul{\Sigma})$
be a family giving a minimal versal deformation of $\mul{\Sigma}$ of
fixed type, and $\mul{\Sigma}_{S_{\mul{\delta}}(\mul{\delta})} \to
S_{\mul{\delta}} \subset \Def(\mul{\Sigma}^{\mul{\delta}})$ a family
giving a minimal versal deformation of $\mul{\Sigma}^{\mul{\delta}}$.
The gluing construction \eqref{glue1} applied to the family
$\mul{\Sigma}_S$ produces a map
\begin{equation} \label{linglue} 
\Def_\Gamma(\mul{\Sigma}) \to \J(\mul{\Sigma}^{\mul{\delta}}), \quad
\zeta \mapsto j^{\mul{\delta}}(\zeta)
\end{equation} 
which maps any deformation of the original curve to the corresponding
deformation of the glued curve.  In other words, any variation of
complex structure on $\mul{\Sigma}$ of fixed type induces a variation
of complex structure on $\mul{\Sigma}^{\mul{\delta}}$.  Similarly, for
any $\xi \in \Omega^0(\mul{\Sigma}, \mul{u}^* TX)$ we obtain an element
$\xi^{\mul{\delta}} \in \Omega^0(\mul{\Sigma}^{\mul{\delta}}, \mul{u}^*
TX)$.  

\begin{proposition} \label{errorprop2}
Suppose that $\mul{u}, \mul{u}^{\mul{\delta}}$ are as above, and
$(\zeta,\xi) \in \Def_\Gamma(u)$.  There is a constant $c$
and an $\eps > 0$ such that if
$ \Vert \mul{\delta} \Vert < \eps, \rho > 1/\eps$, $\Vert \zeta \Vert +
\Vert \mul{\xi} \Vert_{1,p} \leq \eps$, and $ |\delta_k|^2 \rho < \eps$ for
$ k= 1,\ldots, m $ then
$$\Vert \ol{\partial}_{\mul{j}^{\mul{\delta}}(\zeta)}
  \exp_{\mul{u}^{\mul{\delta}}}(\mul{\xi}^{\mul{\delta}})
  \Vert_{0,p,\mul{\delta}}^p \leq c \sum_{k=1}^m (|\delta_k|^{1/2}
  \rho)^{2} .$$
\end{proposition}

\noindent {\em Step 2: Uniformly bounded right inverse}

We wish to show that the map in Proposition \ref{errorprop2} can be
corrected to obtain a holomorphic map.  Define
\begin{equation} \label{cFd}
 \cF_{\mul{u}}^{\mul{\delta}}: 
\Def_\Gamma(\mul{\Sigma})
 \times \Omega^0(\mul{\Sigma}^{\mul{\delta}},\mul{u}^{\mul{\delta},*} TX)
 \to \Omega^{0,1}(\mul{\Sigma}^{\mul{\delta}},\mul{u}^{\mul{\delta},*}
 TX) \end{equation} 
$$ (\zeta,\xi) \mapsto \Psi_{j,\mul{u}^{\mul{\delta}}}(\zeta,\xi)^{-1}
(\olp_{\mul{j}^{\mul{\delta}}(\zeta)}(\exp_{\mul{u}^{\mul{\delta}}}(\xi))) .$$
Here the operator $\Psi_{j,\mul{u}^\delta}$ is as in
\eqref{eq:para}.  Let $\ti{D}_{\mul{u}}^{\mul{\delta}}(\xi)$ be the
associated linear operator, that is, the linearization of \eqref{cFd}
at $\xi$.  This operator naturally extends to a map from Sobolev
$1,p$-completion of the second factor of the domain to the
$0,p$-completion of the codomain.  We denote by
$\ti{D}_{\mul{u}}^{\mul{\delta}} := \ti{D}_{\mul{u}}^{\mul{\delta}}(0)$.
We will construct an approximate inverse
\begin{equation} 
\label{eq:approx}
T_{\mul{\delta}} :
\Omega^{0,1}(\mul{\Sigma}^{\mul{\delta}},\mul{u}^{\mul{\delta},*} TX ) \to
\Def_\Gamma(\mul{\Sigma}) \oplus \Omega^0(\mul{\Sigma}^{\mul{\delta}},\mul{u}^{\mul{\delta},*}
TX)
\end{equation}
to $\ti{D}_{\mul{u}}^{\mul{\delta}}$.  The
construction depends on a carefully chosen cutoff function:

\begin{lemma}\label{careful} \cite[Section 10.3]{ms:jh} 
For any $\delta > 0, \rho > 1$ there exists a function
$\beta_{\rho,\delta} : \R^2 \to [0,1]$ that satisfies
$$ \beta_{\rho,\delta}(z) = \begin{cases} 0 & | z| \leq
  \sqrt{\delta /\rho} \\ 1 & | z| \geq \sqrt{\delta \rho}
\end{cases} $$
and for all $ \xi \in W^{1,p}(B_{\rho |\delta_k|})$ satisfying $
\xi(0) = 0$
\begin{equation} \label{beta1p}
 \Vert (\nabla \beta_{2;\rho,\delta}) \xi \Vert_{0,p} \leq 
c
 \log( \rho^2 )^{-1 + 1/p} \Vert \xi \Vert_{1,p}
 , \quad 
 \Vert \beta_{2;\rho,\delta} \Vert_{1,2} \leq C \log(\rho^{2})^{-1/2}.
 \end{equation}
\end{lemma}

Recall the map $\mul{u}_0$ from Remark \ref{nodalrem}. 

\begin{lemma} \label{rightex}
 For sufficiently small $\mul{\delta}$ there exists a right inverse
$Q_{\mul{u}_0^{\mul{\delta}}}$ of $\ti{D}_{\mul{u}_0^{\mul{\delta}}}$
  with image the $L^2$-perpendicular of the kernel of
  $\ti{D}_{\mul{u}_0^{\mul{\delta}}}$.
\end{lemma} 

\begin{proof} 
Consider the maps defined by parallel transport using the modified
Levi-Civita connection,
$\Pi_{\mul{u}_0^{\mul{\delta}}}^{\mul{u}}: \Omega^0(\mul{\Sigma},
(\mul{u}_0^{\mul{\delta}})^* TX) \to \Omega^0(\mul{\Sigma},
(\mul{u}^{\mul{\delta}})^* TX) $
$ \Psi_{\mul{u}_0^{\mul{\delta}}}^{\mul{u}}: \Omega^{0,1}(\mul{\Sigma},
(\mul{u}_0^{\mul{\delta}})^* TX) \to \Omega^{0,1}(\mul{\Sigma}, (\mul{u}^{\mul{\delta}})^* TX) .$
The operator $ \Psi_{\mul{u}_0^{\mul{\delta}}}^{\mul{u}}
\ti{D}_{\mul{\Sigma},\mul{u}_0^{\mul{\delta}}}
(\Pi_{\mul{u}_0^{\mul{\delta}}}^{\mul{u}})^{-1} $ approaches the
operator $\ti{D}_{\mul{u}}$ as $\mul{\delta} \to 0$, c.f.
\cite[Remark 10.2.2]{ms:jh}.  The statement of the lemma
follows.  \end{proof}

Define the approximate right inverse for
$\ti{D}_{\mul{u}}^{{\mul{\delta}}}$ by composing the right inverse
$Q_{\mul{u}_0^{\mul{\delta}}}$ with a cutoff and extension operator:
$ T_{\mul{\delta}} := P_{\mul{\delta}} Q_{\mul{u}_0^{\mul{\delta}}}
K_{\mul{\delta}} ,$
defined as follows.  The {\em cutoff operator}
$$ K_{\mul{\delta}}: \Omega^{0,1}(\mul{u}^{\mul{\delta}}* TX)_{0,p,\mul{\delta}} \to
\Omega^{0,1}(\mul{u}_0^{\mul{\delta},*} TX)_{0,p} $$
is defined by
$$ (K_{\mul{\delta}}(\eta))(z) =
\begin{cases}  \eta(z) & z \notin\bigcup_{k,\pm} B_{| \delta_k(z)|^{1/2}} (w_k^\pm)  \\
                0 & \text{otherwise} \end{cases} .$$
We have 
$ \Vert  K_{\mul{\delta}} \eta \Vert_{0,p}  
\leq \Vert \eta \Vert_{0,p,\mul{\delta}} $
by definition of the $0,p,\mul{\delta}$ norm.  The {\em extension operator}
$$ P_{\mul{\delta}}: \Def_\Gamma(\mul{\Sigma}) \oplus
  \Omega^0(\mul{\Sigma},\mul{u}_0^{\mul{\delta}}* TX)_{1,p,\mul{\delta}} \to
  \Def(\mul{\Sigma}^{\mul{\delta}}) \oplus \Omega^0(\mul{\Sigma}^{\mul{\delta}},\mul{u}^{\mul{\delta},*}
  TX)_{1,p,\mul{\delta}} $$
is defined as follows.  For each component $\Sigma_i$ let $\Sigma_i^*$
denote the complements of small balls around the nodes
$$ \Sigma_i^* = \Sigma_i - \bigcup_{l, w_l^\pm \in 
\Sigma_i}  B_{ | \delta_l(z)|^{1/2}/\rho} (w_l^\pm)  $$
and the inclusion $\pi_i : \Sigma_i^* \to \mul{\Sigma}^{\mul{\delta}} $
induces a map $ \pi_{i,*} : \Omega^0( \Sigma_i^*,u_i^* TX) \to
\Omega^0(\mul{\Sigma}^{\mul{\delta}}, \mul{u}^{\mul{\delta},*} TX )_{0,p} .$
Define
$$ P_{\mul{\delta}}(\zeta,\xi) =
(\zeta^{\mul{\delta}},\xi^{\mul{\delta}}) $$
where $ \zeta^{\mul{\delta}}$ is the image of $\zeta$ under the gluing
map \eqref{linglue} and $\xi^{\mul{\delta}}$ is obtained by patching
together the sections $\xi$; on the gluing region arising from gluing
the $k$-th node $w_k$ the section $\xi^{\mul{\delta}}$ is given by the
sum
$$ \pi_{i^+(k),*} 
\beta_{\rho,\delta_k} (\xi_{i^+(k)} -
 \xi(w_k)) \\ + \pi_{i^-(k),*}  \beta_{\rho,\delta_k} (\xi_{i^-(k)} -
 \xi(w_k)) + \xi(w_k).
$$
Fix a metric $\Vert \cdot \Vert$ on the finite-dimensional space
$\Def_\Gamma(\mul{\Sigma})$ and define
$$\Vert (\zeta,\xi) \Vert_{1,p,\mul{\delta}} = \left( \Vert \zeta \Vert^p
+ \Vert \xi \Vert_{1,p,\mul{\delta}}^p \right)^{1/p} .$$

\begin{proposition} \label{right}
Let $\mul{u} : \mul{\Sigma} \to X$ be a stable map.  There exist
constants $c,C > 0 $ such that if $ \Vert \mul{\delta} \Vert < c $ then
the approximate inverse $ T_{\mul{\delta}}$ of  \eqref{eq:approx}
satisfies
$$ \Vert (\ti{D}_{\mul{u}}^{\mul{\delta}} T_{\mul{\delta}} - I) \eta
\Vert_{0,p,\mul{\delta}} \leq \hh \Vert \eta \Vert_{0,p,\mul{\delta}}, 
\quad  \Vert T_{\mul{\delta}} \Vert < C .$$ \end{proposition}

\begin{proof}  By construction $T_{\mul{\delta}}$ is an exact right
inverse for $\ti{D}_{\mul{u}}^{\mul{\delta}}$ away from gluing region.  In the
gluing region the variation of complex structure on the curve vanishes
and $D_{\mul{u}}^{\mul{\delta}} = D_{x_k}$, the standard Cauchy-Riemann operator
with values in $T_{x_k}X$.  So
\begin{eqnarray*} \ti{D}_{\mul{u}}^{\mul{\delta}} T_{\mul{\delta}} \eta - \eta &=& \sum D_{x_k}
\beta_{2;\rho,\delta}(z) ( \xi_{i^\pm(k)}(z) - \xi_{i^\pm(k)}(w_k)) \\
&=& \sum (D_{x_k} \beta_{2;\rho,\delta}(z) ) \xi_{i^\pm(k)}(z) 
\end{eqnarray*}
since $K_{\mul{\delta}} \eta = 0 $ on $B_{| \delta_k|^{1/2}}(0)$ in the
components adjacent to the node.  Since $p > 2$, the
$0,p,\mul{\delta}$-norm of the right hand side is controlled by the
ordinary $L^p$ norm.  By \eqref{beta1p} we have
$$  \Vert \ti{D}_{\mul{u}}^{\mul{\delta}} T_{\mul{\delta}} \eta - \eta \Vert_{0,p,\mul{\delta}} \leq \sum_k c
| \log( \rho )|^{2/p-2} \Vert \xi_{i^\pm(k)}
- \xi(w_k) \Vert_{1,p}
.$$
The last factor is bounded by $ \Vert K_{\mul{\delta}} \eta \Vert_{0,p}$, by
the uniform bound on $Q_{\mul{\delta}}$, and hence $\Vert \eta
\Vert_{0,p,\mul{\delta}}$, by the uniform bound on $K_{\mul{\delta}}$.
\end{proof}  

Define a right inverse $Q_{\mul{\delta}}$ to
$\ti{D}_{\mul{u}}^{\mul{\delta}}$ by the formula
$$ Q_{\mul{\delta}} = T_{\mul{\delta}} (\ti{D}_{\mul{u}}^{\mul{\delta}} T_{\mul{\delta}})^{-1} = \sum_{k \ge 0}
T_{\mul{\delta}} (\ti{D}_{\mul{u}}^{\mul{\delta}} T_{\mul{\delta}} - I )^k .$$
The uniform bound on $T_{\mul{\delta}}$ from Lemma \ref{right} implies
a uniform bound on $Q_{\mul{\delta}}$.

\noindent {\em Step 3: Uniform quadratic estimate}

\begin{proposition}
\label{quadratic}  
Let $\mul{u}: \mul{\Sigma} \to X$ be a stable map, $\mul{\delta}$ a
collection of gluing parameters, and $\mul{u}^{\mul{\delta}}:
\mul{\Sigma}^{\mul{\delta}} \to X$ the approximate solution defined
above.  For every constant $c > 0$ there exist constants $c_0,\delta_0
> 0$ such that if $ \mul{u} \in \Map(\mul{\Sigma},X)_{1,p}, \xi \in
\Omega^0(\mul{\Sigma}^{\mul{\delta}},\mul{u}^{\mul{\delta},*}TX)_{1,p}$
$$ \Vert \d \mul{u} \Vert_{0,p} \leq c_0, \ \ \ \Vert \xi
\Vert_{L_\infty} \leq c_0 \ \ \ \ \Vert \zeta \Vert \leq c_0,
\ \ \ \ \Vert \mul{\delta} \Vert < \delta_0
$$
then
$$ \Vert (D \cF_{\mul{u}}^{\mul{\delta}}(\zeta,\xi, \zeta_1,\xi_1) -
\ti{D}_{\mul{u}}^{{\mul{\delta}}})(\zeta_1,\xi_1)
\Vert_{0,p,\mul{\delta}} \leq c \Vert \zeta, \xi \Vert_{1,p} ( \Vert
\zeta_1, \xi_1 \Vert_{1,p} .$$
\end{proposition} 

\noindent Here $D \cF_{\mul{u}}^{\mul{\delta}}(\zeta,\xi,
\zeta_1,\xi_1)$ denotes the derivative evaluated at $\zeta,\xi$,
applied to $\zeta_1,\xi_1$.  We use similar notation throughout the
discussion.  The proof uses uniform estimates on Sobolev embedding:

\begin{lemma} \label{unifsob1} There exists
a constant $c > 0$ independent of $\delta$ such that the embedding
$$ 
\Omega^0(\mul{\Sigma}^{\mul{\delta}}, \mul{u}^{\mul{\delta},*} 
T X)_{{1,p},\mul{\delta}} 
\to 
\Omega^0(\mul{\Sigma}^{\mul{\delta}}, \mul{u}^{\mul{\delta},*} 
T X)_{0,\infty} $$
has norm less than $c$.
\end{lemma}  

\begin{proof}  One writes the Sobolev norms as a contribution from 
each component of the curve $\mul{\Sigma}$.  Then on each piece,
the metric near the boundary is uniformly comparable with
the flat metric.  The claim then follows from 
\cite[Chapter 4]{ad:so} which shows that the constants 
in the Sobolev embeddings depend only on the dimensions of the cone in
the cone condition.
\end{proof}

\begin{proof}[Proof of Proposition]  For simplicity we assume a single gluing parameter $\delta$.
Let
$ \Psi_{\mul{u}}^{{\delta},x}(\zeta,\xi): \Lambda^{0,1} T_z^*
\mul{\Sigma}^{\delta} \otimes T_x X \to \Lambda^{0,1}_{\mul{j}^{\mul{\delta}}(\zeta)} T_z^*
\mul{\Sigma} \otimes T_{\exp_x(\xi)} X $
denote pointwise parallel transport as in \eqref{parallel} using the
parallel transport using the modified Levi-Civita connection, and
projecting onto the $0,1$-part of the form defined using the complex
structure $\mul{j}^{\mul{\delta}}(\zeta)$ obtained from gluing
$\mul{j}(\zeta)$, see \eqref{linglue}.  Let
$$ \Theta_{\mul{u}}^{{\delta},x}(\zeta,\xi,\zeta_1,\xi_1;\eta) =
  \ti{\nabla}_t \Psi_{\mul{u}}^{{\delta},x} ( \zeta + t \zeta_1, \xi +
    t\xi_1) \eta .$$
For $\xi, \eta$ sufficiently small there exists a constant $c$ such
that
\begin{equation} \label{Psiest}
 | \Theta_{\mul{u}}^{\delta,x}(\zeta,\xi,\zeta_1,\xi_1; \eta )| \leq c \Vert
\xi,\zeta \Vert \Vert \xi_1,\zeta_1 \Vert \Vert \eta \Vert
\end{equation}
where the norms on the right-hand side are any norms on the finite
dimensional vector spaces $T_{\mul{\Sigma}} M_{g,n,\Gamma}$ and $T_x
X$.  This estimate is uniform in $\delta$, for the variation in
complex structure vanishes in a neighborhood of the nodes.
Differentiate the equation
$ \Psi_{\mul{u}}^{\delta}(\zeta,\xi) \cF_{\mul{u}}^{\delta}(\zeta,\xi)
= \olp_{\mul{j}^\delta(\zeta)}(\exp_{\mul{u}^\delta}(\xi))) $
with respect to $(\zeta_1,\xi_1)$ to obtain
\begin{multline}
 \Theta_{\mul{u}^\delta}(\zeta,\xi,\zeta_1,\xi_1,
 \cF_{\mul{u}}^{\delta}(\zeta,\xi) ) +
 \Psi_{\mul{u}}^{\delta}(\zeta,\xi)( D
 \cF_{\mul{u}}^{\delta}(\zeta,\xi,\zeta_1,\xi_1)) =
 \\ \ti{D}_{\mul{u}}^{\delta}(\xi, Dj^\delta
 (\zeta,\zeta_1),D\exp_{\mul{u}^\delta} (\xi,\xi_1)) .\end{multline}
Using the pointwise inequality 
$$ | \cF_{\mul{u}}^{\delta}(\zeta,\xi) | < c | \d
{\exp_{\mul{u}^\delta}(\xi)} | < c ( | \d \mul{u}^\delta | + | \nabla
\xi | )
$$
for $\zeta,\xi$ sufficiently small, the estimate \eqref{Psiest} on
$\Phi$ produces a pointwise estimate
$$ | (\Psi_{\mul{u}}^{\delta})^{-1}(\xi)
\Theta_{\mul{u}}^{\delta}(\zeta,\xi,\zeta_1,\xi_1,\cF_{\mul{u}}^{\delta}(\zeta,\xi))
| \leq c (| \d \mul{u}^\delta | + | \nabla \xi |) \, | ( \xi,\zeta ) |
\, | (\xi_1,\zeta_1) | .$$
Hence
\begin{multline} \Vert \Psi_{\mul{u}^\delta}^{-1}(\xi)
\Theta_{\mul{u}}^{\delta}(\zeta,\xi,\zeta_1,\xi_1,\cF^\delta_{\mul{u}}(\zeta,\xi))
\Vert_{0,p} \\ \leq c ( 1+ \Vert \d \mul{u}^\delta \Vert_{0,p} + \Vert \nabla \xi
\Vert_{0,p} ) \Vert (\xi,\zeta)
\Vert_{0,\infty} \Vert (\xi_1,\zeta_1) \Vert_{0,\infty}
.\end{multline}
It follows that 
\begin{equation} \label{firstclaim}
  \Vert \Psi_{\mul{u}}^{\delta}(\xi)^{-1}
\Theta_{\mul{u}}^{\delta}(\zeta,\xi,\zeta_1,\xi_1,\cF_{\mul{u}}^{\delta}(\zeta,\xi)) \Vert_{0,p}
\leq c \Vert (\xi,\zeta) \Vert_{1,p} \Vert (\xi_1,\zeta_1) \Vert_{1,p}
\end{equation}
since the $W^{1,p}$ norm controls the $L^\infty$ norm by Lemma \ref{unifsob1}. 

We next show that there exists a constant $c > 0$ such that uniformly
in $\delta$,
\begin{equation} \label{secondclaim} \Vert \Psi_{\mul{u}^\delta}(\xi)^{-1}
\ti{D}_{\mul{u}}^\delta(\xi, D\mul{j}^{\mul{\delta}}(\zeta,\zeta_1),D
\exp_{{\mul{u}^\delta}} (\xi,\xi_1)) - \ti{D}_{\mul{u}}^{\delta}
(\zeta_1,\xi_1) \Vert_{0,p}
\leq c \Vert \zeta,\xi \Vert_{1,p} \Vert \zeta_1,\xi_1 \Vert_{1,p}
.\end{equation}
Indeed differentiate \eqref{cFd} to obtain
\begin{multline}
 \ti{D}_{\mul{u}}^{\delta}(\xi, D \mul{j}^\delta(\zeta,\zeta_1),D
 \exp_{{\mul{u}^\delta}}(\xi,\xi_1)) = \nabla^{0,1}_{j^\delta(\zeta)}
 D\exp_{{\mul{u}^\delta}}(\xi,\xi_1) - \\ \hh \pi^{0,1}_{j^\delta(\zeta)}
 J_{\exp_{\mul{u}^\delta}(\xi)} \d \exp_{\mul{u}^\delta}(\xi)
 D\mul{j}^\delta(\zeta,\zeta_1) $$ - J_{\exp_{\mul{u}^\delta}(\xi)} (\nabla_{D
   \exp_{{\mul{u}}}(\xi,\xi_1)} J_{\exp_{\mul{u}^\delta}(\xi)} )
 \partial \exp_{\mul{u}^\delta}(\xi) . \end{multline}
Hence 
$$ \ti{D}_{\mul{u}}^\delta(\xi, D\mul{j}^\delta(\zeta,\zeta_1),D
\exp_{{\mul{u}^\delta}}(\xi,\xi_1)) 
- \Psi_{\mul{u}^\delta}(\xi) \ti{D}_{\mul{u}}^{\delta}(\zeta_1,\xi_1) 
= \Pi_1 + \Pi_2 + \Pi_3 $$
where the three terms $\Pi_1,\Pi_2,\Pi_3$ are
$$ \Pi_1 = 
\nabla^{0,1}_{\mul{j}^\delta(\zeta)} D\exp_{{\mul{u}^\delta}}(\xi,\xi_1)
- 
\Psi_{\mul{u}^\delta}(\xi) \nabla^{0,1}_{\mul{j}^\delta(0)} \xi_1$$
$$ \Pi_2 = - \hh \pi^{0,1}_{\mul{j}^\delta(\zeta)} J_{\exp_{\mul{u}}(\xi)} \d
\exp_{\mul{u}^\delta}(\xi) D\mul{j}^\delta(\zeta,\zeta_1) + 
\Psi_{\mul{u}^\delta}(\xi) \pi^{0,1}_{\mul{j}^\delta(0)} \hh J_{\mul{u}^\delta} \d
\mul{u}^\delta D\mul{j}^\delta(0,\zeta_1)
$$ 
$$ \Pi_3 = - \hh J_{\exp_{\mul{u}}(\xi)} (\nabla^{0,1}_{D
  \exp_{{\mul{u}}^\delta}(\xi,\xi_1)} J_{\exp_{\mul{u}^\delta}(\xi)} ) \partial_{\mul{j}^\delta(\zeta)}
\exp_{\mul{u}^\delta}(\xi) + \hh \Psi_{\mul{u}^\delta}(\xi) J_{\mul{u}^\delta}
(\nabla_{\xi_1} J_{\mul{u}} ) \partial_{\mul{j}^\delta(0)} \mul{u}^\delta
.$$
The first difference has norm given by 
\begin{multline}   | \pi_{\mul{j}^\delta(\zeta)}^{0,1} ( \nabla ( D \exp_{{\mul{u}^\delta}}(\xi,\xi_1)) 
- \Psi_{\mul{u}^\delta}(\xi) \nabla \xi_1 ) | \\ \leq |
\pi_{\mul{j}^\delta(\zeta)}^{0,1} ( \nabla ( D \exp_{{\mul{u}^\delta}}(\xi,\xi_1)) -
D\exp_{{\mul{u}^\delta}}(\xi,\nabla \xi_1) )| + |\pi^{0,1}_{\mul{j}^\delta(\zeta)} (
D\exp_{{\mul{u}^\delta}}(\xi,\nabla \xi_1) - \Psi_{\mul{u}^\delta}(\xi)
\nabla \xi_1 ) | \\ \leq c | \nabla \xi | | \xi_1| + c ( | \zeta| + |
\xi |) |\nabla \xi_1| + c | \d \mul{u}^\delta | | \xi| | \xi_1|
.\end{multline}
We write for the second difference
\begin{multline}
  | \pi^{0,1}_{\mul{j}^\delta(\zeta)}( J_{\exp_{\mul{u}^\delta}(\xi)} \d \exp_{\mul{u}^\delta}(\xi) D \mul{j}^\delta(\zeta,\zeta_1) -
\Psi_{\mul{u}^\delta}(\xi) J_{\mul{u}^\delta} \d \mul{u}^\delta D\mul{j}^\delta(0,\zeta_1))|  \\
\leq 
 | \pi^{0,1}_{\mul{j}^\delta(\zeta)} ( J_{\exp_{\mul{u}^\delta}(\xi)} \d \exp_{\mul{u}^\delta}(\xi) D\mul{j}^\delta(\zeta,\zeta_1) -
  J_{\exp_{\mul{u}^\delta}(\xi)}  \Psi_{\mul{u}^\delta}(\xi)  \d\mul{u}^\delta D\mul{j}^\delta(0,\zeta_1))|  
+ \\
|\pi^{0,1}_{\mul{j}^\delta(\zeta)} ( J_{\exp_{\mul{u}^\delta}(\xi)}  \Psi_{\mul{u}^\delta}(\xi)  \d\mul{u}^\delta D\mul{j}^\delta(0,\zeta_1)   - 
\Psi_{\mul{u}^\delta}(\xi) J_{\mul{u}} \d \mul{u}^\delta D\mul{j}^\delta(0,\zeta_1))|  \\
 \leq c ( |\zeta| + |\xi|  + | \d \mul{u}^\delta| +  | \nabla \xi |)  | |\zeta_1|  .
\end{multline}
The third term can be estimated pointwise by
$$ | J_{\exp_{\mul{u}^\delta}(\xi)} 
(\nabla_{D \exp_{\mul{u}^\delta}(\xi,\xi_1)}
J_{\exp_{\mul{u}^\delta}(\xi)} ) \partial_{\mul{j}^\delta(\zeta)} \exp_{\mul{u}^\delta}(\xi) -
\Psi_{\mul{u}^\delta}(\xi) J_{\mul{u}^\delta}
(\nabla_{\xi_1} J_{\mul{u}^\delta} ) \partial_{\mul{j}^\delta(0)}
\mul{u}^\delta| $$
$$ \leq | J_{\exp_{\mul{u}^\delta}(\xi)} 
(\nabla_{D \exp_{{\mul{u}}}(\xi,\xi_1)} J_{\exp_{\mul{u}^\delta}(\xi)} )
\partial_{\mul{j}^\delta(\zeta)} \exp_{\mul{u}^\delta}(\xi) - 
J_{\exp_{\mul{u}^\delta}(\xi)} (\nabla_{D
  \exp_{\mul{u}^\delta}(\xi,\xi_1)} J_{\exp_{\mul{u}^\delta}(\xi)} )
\partial_{\mul{j}^\delta(0)} \exp_{\mul{u}^\delta}(\xi)| $$
$$ + 
| J_{\exp_{\mul{u}^\delta}(\xi)} 
(\nabla_{D \exp_{{\mul{u}}}(\xi,\xi_1)} J_{\exp_{\mul{u}^\delta}(\xi)} 
) \partial_{\mul{j}^\delta(0)} \exp_{\mul{u}^\delta}(\xi)  
- 
\Psi_{\mul{u}^\delta}(\xi)  J_{\mul{u}^\delta} 
 (\nabla_{\xi_1} J_{\mul{u}^\delta} 
) \partial_{\mul{j}^\delta(0)} \mul{u}^\delta|   $$
$$ \leq c | \zeta | (| \d \mul{u}^\delta| + |\nabla \xi |) |\xi_1| + c ( |
\d \mul{u}^\delta | + | \nabla \xi |) ( | \xi_1 |) $$
for $\xi$ sufficiently small.  Combining these estimates and
integrating, using the $0,p,\delta$-norms on $\d \mul{u}$, $\nabla \xi,
\nabla \xi_1$ and the $L^\infty$ norms on the other factors and Lemma
\ref{unifsob1}, completes the proof.
\end{proof} 

\noindent {\em Step 4: Implicit Function Theorem}

For any $(\zeta_0,\xi_0) \in \Def_\Gamma (\mul{u})$ we denote by
$\zeta_0^{\mul{\delta}}$ the deformation of $\mul{\Sigma}^{\mul{\delta}}$
defined in \eqref{linglue} and by $\xi_0^{\mul{\delta}}$ the section of
$\mul{u}^{\mul{\delta},*} TX$ defined as in \eqref{preglued}.

\begin{theorem} \label{gluingstable}
Let $\mul{u}:\mul{\Sigma} \to X$ be a stable map. 
 There exist constants
$\eps_0,\eps_1 > 0$ such that for any $(\zeta_0,\xi_0,\mul{\delta}) \in
\ker \ti{D}_{\mul{u}} \times \R^m$ of norm at most $\eps_0$, there is a
unique
$(\zeta_1,\xi_1) = (\ti{D}_{\mul{u}}^{\mul{\delta}})^* \eta_1 $ 
of norm at most $\eps_1$ such that the map
$\exp_{\mul{u}^{\mul{\delta}}}(\xi_0^{\mul{\delta}} + \xi)$ is
$\mul{j}^\delta(\zeta_0 + \zeta_1)$-holomorphic, and depends smoothly on
$\zeta_0,\xi_0$.
\end{theorem}  

\begin{proof}  The first claim is an application of the quantitative version of the implicit
function theorem (see for example \cite[Appendix A.3]{ms:jh}) using
the uniform error bound from Proposition \ref{errorprop}, uniformly
bounded right inverse from Proposition \ref{right}, and uniform
quadratic estimate from Proposition \ref{quadratic}.  
\end{proof}

\noindent {\em Step 5: Rigidification} 

In the previous step we have constructed a family of stable maps which
we will show eventually gives rise to a parametrization of all nearby
stable maps.  A more natural way of parametrizing nearby stable maps
involves examining the intersections with a family of codimension two
submanifolds.  For example, this construction of charts is that given
in the algebraic geometry approach of Fulton-Pandharipande
\cite{fu:st}.  In order to carry this out in the symplectic approach,
we study the differentiability of the evaluation maps.  Let $u_S:
\Sigma_S \to X$ over a parameter space $S \subset \Def(\mul{u})$ be
the family of maps defined in the previous step.  The following is
similar to \cite[Lemma A1.59]{fooo}.

\begin{theorem} \label{diffev} If $ \mul{u}_S$
is constructed using the exponential gluing profile $\varphi$ and $U
\subset \mul{\Sigma}$ is an open neighborhood of the nodes then the
map $(z,s) \mapsto \mul{u}_s(z)$ is differentiable on a neighborhood
of $(\mul{\Sigma} - U ) \times \{ 0 \}$.
\end{theorem}  

\begin{proof} 
 For simplicity, we assume that there is a single gluing parameter
 $\delta$.  Differentiability for $\delta$ is studied in
 McDuff-Salamon \cite[Section 10.6]{ms:jh}.  The discussion in our
 case is somewhat easier, because we use a fixed right inverse in the
 gluing construction.  Given $(\zeta,\xi) \in \Def_\Gamma(\mul{\Sigma})
 \times \Omega^0( \mul{\Sigma}^\delta, (\mul{u}^\delta)^* TX)$, we
 constructed a unique correction $(\zeta_1,\xi_1)$ in the image of the
 right inverse such that
$ \olp_{\mul{j}^\delta(\zeta_0 + \zeta_1)} \exp_{\mul{u}^\delta}(\xi_0 +
 \xi_1) = 0 .$
For $\delta$ fixed, $(\zeta_1,\xi_1)$ depends smoothly on
$(\zeta_0,\xi_0)$, by the implicit function theorem.  Hence the
evaluation at $z \in \mul{\Sigma} - U$ also depends smoothly on
$(\zeta_0,\xi_0)$.

The computation of the derivative with respect to the gluing parameter
is complicated by the fact that for each $\delta$ a different implicit
function theorem is applied to obtain the correction.  Let
$\ti{D}_\delta = D \cF_{\mul{u}^{\delta}}$.  Differentiating the
equation $\cF_{\mul{u}^{\delta}}(\zeta_0^\delta + \zeta_1,\xi_0^\delta
+ \xi_1) = 0 $ with respect to $\delta$ gives
$$ \ti{D}_\delta \left( \ddd \zeta_1, D\exp_{\mul{u}^\delta} (\xi_0^\delta +
\xi_1 , 0 , \ddd \xi_1 ) \right) = - \ti{D}_\delta \left( 0 ,
D\exp_{\mul{u}^\delta} (\xi_0^\delta + \xi_1 , \ddd \mul{u}^\delta , 0 )
\right) .$$
From \eqref{preglued} we
have in the gluing region, 
\begin{eqnarray*}
 \ddd \olp \mul{u}^\delta 
 &=& \ddd \olp \exp_{x} ( \left(  
\alpha( \rho \varphi^{-1/2} |z|) \xi(z)) \right) \\
 &=& 
D\exp_{x} \left(
\alpha( \rho \varphi^{-1/2}  |z| ) \xi(z),
\alpha'( \rho \varphi^{-1/2} |z| )|z| \xi(z)
 \ddd \varphi^{-1/2} \rho \right)^{0,1}  .\end{eqnarray*}
Hence there exists a constant $C$ depending on $\rho,\alpha$ but not
on $\delta$ such that
\begin{equation} \label{point1}
 \left| \ddd \olp \mul{u}^\delta \right| \leq C \left|
\ddd \varphi^{-1/2}  \right| .\end{equation}
Now $\d\varphi$ is given by
$$ \d ( e^{1/\delta} - e)^{-1/2} = (1/2) (e^{1/\delta} - e)^{-3/2} e^{1/\delta}
\delta^{-2} \d \delta = (1/2) ( e^{1/3\delta} - e^{-2/3\delta + 1})^{-3/2} \delta^{-2} \d \delta .$$
For $\delta$  small, this  is less than  $\hh e^{-1/2\delta} \delta^{-2}$.   Integrating and
using the  pointwise  estimate  \eqref{point1}  we  obtain  for  some
constant $C > 0$,
$$ \left\Vert \ddd \olp \mul{u}^\delta_s \right\Vert_{0,p} \leq C
e^{-1/2\delta} \delta^{-2} \leq Ce ^{-1 / \delta}$$
for sufficiently small $\delta$.  Now the uniform quadratic estimates
imply that $\ti{D}_\delta = D \cF_{\mul{u}^\delta}(\zeta,\xi)$ is
uniformly bounded from below on the right inverse of
$\ti{D}_{\mul{u}}^{\delta} = D \cF_{\mul{u}}^{\delta}(0,0)$, for
$(\zeta_0,\xi_0)$ sufficiently small.  It follows that
$$ \left\Vert ( \ddd \zeta_1, \ddd \xi_1) \right\Vert_{1,p} \leq C
e^{-1/\delta} $$
for $\zeta_0,\xi_0,\delta$ sufficiently small as well. Hence the same is
true for the evaluation $\ddd \xi_1(z)$ for $z \in \mul{\Sigma} - U$.
In particular,
$\lim_{\delta \to 0 } (\partial_\delta
\exp_{\mul{u}^\delta}(\xi_0^\delta + \xi_1))(z) = 0.  $
It follows that the differential of the evaluation map has a
continuous limit as ${\delta} \to 0$, which completes the proof of the
Theorem.
\end{proof}

Using the evaluation maps in the previous step, we construct
embeddings of the families constructed above into suitable moduli
spaces of stable marked curves, given by adding additional marked
points which map to fixed submanifolds in $X$.  A codimension two
submanifold $Y \subset X$ is {\em transverse} to $\mul{u}: \mul{\Sigma}
\to X$ if $\mul{u}$ meets $Y$ transversally in a single point
$\mul{u}(z)$.

\begin{definition}  Let $\mul{u}: \mul{\Sigma} \to X$ be a stable map.  
Given any family $\mul{Y} = (Y_1,\ldots, Y_\ell)$ of codimension two
submanifolds transverse to $\mul{u}$ and a family
$\mul{\Sigma}_S,\mul{u}_S,\mul{z}_S$ with parameter space $S$ of an
$n$-marked stable map $(\mul{\Sigma},\mul{u},\mul{z})$, the {\em
  rigidified family} of $n+\ell$-marked nodal surfaces is defined by
\begin{equation} \label{rigidify}  
\mul{\Sigma}_S^{\mul{Y},\mul{u}} := (\mul{\Sigma}_S,
(z_{1,S}, \ldots , z_{n+\ell,S})) \to S,
\ \ \ \ \mul{u}_s(z_{n+i,s}) \in Y_i .\end{equation}
\end{definition} 

\begin{proposition}\label{diffprop}  Let $\mul{u}_S$ be a family of stable
maps over a parameter space $S \subset \Def(\mul{u})$ given by the
gluing construction using a gluing profile $\varphi$ and system of
coordinates $\mul{\kappa}$.  Suppose that the evaluation map $\ev:
(\mul{\Sigma} - U) \times S \to X$ is $C^1$, and that the rigidified
family has stable underlying curves.  Then the rigidified family of
curves $\mul{\Sigma}_S^{\mul{Y},\mul{u}}$ is $C^1$ with respect to the
gluing profile and local coordinates, that is, the map
$ S \mapsto \ol{M}_{g,n+l}, \quad s \mapsto
\mul{\Sigma}_s^{\mul{Y},\mul{u}} $
is $C^1$ with respect to the smooth structure defined by $\varphi,
\mul{\kappa}$.
\end{proposition}  

\begin{proof}  By the implicit function theorem for $C^1$ maps 
and differentiability of evaluation maps from the previous subsection.
\end{proof} 

\begin{definition} \label{compat} 
Let $\mul{Y},\mul{u}$ be as above.  The pair $(\mul{Y},\mul{u})$ is {\em
  compatible} if
\begin{enumerate}
\item each $Y_j$ intersects $\mul{u}$ transversally in a single point
  $z_j \in \mul{\Sigma}$;
\item if $\xi \in \ker(\ti{D}_{\mul{u}})$ satisfies $\xi(z_{n+j}) \in
  T_{\mul{u}(z_{n+j})} Y_j$ for $j =1,\ldots, l$ then $\xi = 0$;
\item the curve $\mul{\Sigma}$ marked with the additional points
  $z_{n+1},\ldots, z_{n+\ell}$ is stable;
\item if some automorphism of $(\mul{\Sigma},\mul{u})$ maps 
$z_i$ to $z_j$ then $Y_i$ is equal to $Y_j$.  
\end{enumerate} 
\end{definition} 

The second condition says that there are no infinitesimal deformations
which do not change the positions of the extra markings.  

\begin{proposition} \label{versalrigid} 
Let $\mul{u}$ be a parametrized regular stable map, and $\mul{u}_S$ the
stratified-smooth universal deformation constructed in Theorem
\ref{gluingstable} with base $S \subset \Def(\mul{u})$.  There exists a
collection $\mul{Y}$ compatible with $\mul{u}$.  Furthermore if the
evaluation map is $C^1$ as in Proposition \ref{diffprop} then $
\mul{\Sigma}^{\mul{Y},\mul{u}}_S$ defines an $C^1$-immersion of $S$ into
$\Def(\mul{\Sigma}^{\mul{Y},\mul{u}})$.
\end{proposition}

\begin{proof}  
First we show the existence of a compatible collection.  Given a
regular stable map $(\mul{\Sigma},\mul{z} =
(z_1,\ldots,z_n),\mul{u}:\mul{\Sigma} \to X)$, choose $Y_1,\ldots,Y_k$
transverse $\mul{u}$ on the unstable components of $\mul{\Sigma}$, so
that $\mul{\Sigma}_1 = (\mul{\Sigma}, (z_1,\ldots,z_{n+k}))$ is a
stable curve.  Let $\mul{\Sigma}_{S_1,1} \to S_1$ denote a universal
deformation of $\mul{\Sigma}_1$.  By universality, the family
$\mul{\Sigma}_S^{\mul{Y},\mul{u}}$ is induced by a map $\psi : S \to
S_1$.  We successively add marked points until $\psi$ is an immersion:
Suppose that $\psi$ is not an immersion.  Then we may choose an
additional marked point $z_{n+k+1} \in \mul{\Sigma}$ such that $\d
\ev_{n+k+1}$ is non-trivial on $\ker D \psi$.  Since $\mul{u}$ is
holomorphic, $\d \mul{u}(z_{n+k+1})$ is rank two at $z_{n+k+1}$.  Let
$Y_{n+k+1} \subset X$ be a codimension two submanifold containing
$\mul{u}(z_{n+k+1})$ such that $\mul{u}$ is transverse to $Y_{n+k+1}$
at $z_{n+k+1}$, and $Y_{n+k+1}$ is transversal to $\ev_{n+k+1} $ at
$\mul{\Sigma},\mul{u}$.  Suppose $z_{n+n'+1}$ has orbit $z_{n+k+1},
z_{n+k+2},\ldots, z_{n+l}$ under the group $\Aut(\mul{u})$.  Repeating
the same submanifold for each marking related by automorphisms gives
a collection invariant under the action of automorphisms.  The map
$\psi_1$ for the new family has property that the dimension of
$\ker(D\psi_1)$ has dimension at least two less than that of
$\ker(D\psi)$.  It follows that the procedure terminates after adding
a finite number of markings.  The last claim follows from the second
condition in Definition \ref{compat}.  \end{proof}

\noindent {\em Step 6: Surjectivity}

In this section, we show that the family constructed above contains a
Gromov neighborhood of the central fiber.  First we show:

\begin{proposition}  \label{surject} There exists a constant $\eps > 0$ 
such that any stable map $(\mul{\Sigma}_1,\mul{u}_1)$ with complex
structure on $\mul{\Sigma}_1$ given by $j^\delta(\zeta)$ for some
$\zeta \in \Def(\mul{\Sigma}^{\mul{\delta}})$, and $ \mul{u}_1 :=
\exp_{\mul{u}^{\mul{\delta}}}(\mul{\xi}) $ with
$ \Vert \zeta \Vert^2 + \Vert \xi \Vert^2_{1,p,\mul{\delta}} < \eps $
is of the form in Theorem \ref{gluingstable} for some $(\xi_1,\zeta_1)
\in \Im (\ti{D}_{\mul{u}}^{{\mul{\delta}}})^*$.  
\end{proposition} 

\begin{proof}   Compare with \cite[Section 10.7.3]{ms:jh}.  
 Let $(\zeta,\mul{\xi})$ be as in the statement of the Proposition.  We
 claim that $(\zeta,\mul{\xi}) =
 (\zeta_0^{\mul{\delta}},\mul{\xi}_0^{\mul{\delta}}) +
 (\zeta_1,\mul{\xi}_1)$ for some $(\zeta_0,\mul{\xi}_0) \in
 \ker(\ti{D}_{\mul{u}})$ and $(\zeta_1,\mul{\xi}_1) \in \Im
 ((D_{\mul{u}}^{{\mul{\delta}}})^*)$ with small norm.  It then follows
 by the implicit function theorem that $(\zeta_1,\mul{\xi}_1)$ is the
 solution given in Theorem \ref{gluingstable}.  Now $\zeta =
 \zeta_0^{\mul{\delta}}$ for some $\zeta_0 \in
 \Def_\Gamma(\mul{\Sigma})$ and gluing parameters $\mul{\delta}$,
 because $\Def(\mul{\Sigma}^{\mul{\delta}})$ is the direct sum of the
 image of $\Def_\Gamma(\mul{\Sigma})$ and $\C^m$.  By the gluing
 theorem for indices (see e.g. \cite[Theorem 2.4.1]{orient}), the
 image of $\Def_\Gamma(\mul{u})$ under the gluing map projects
 isomorphically onto $\ker(\ti{D}_{\mul{u}}^{\mul{\delta}})$ for
 $\mul{\delta}$ sufficiently small, and so $\Def_\Gamma(\mul{u})$ is
 transverse to $\Im (\ti{D}_{\mul{u}}^{\mul{\delta}})^*$, for
 $\mul{\delta}$ sufficiently small.  The claim then follows from the
 inverse function theorem.
\end{proof}

Given a regular stable $\mul{u}$ with stable domain, consider the
family of $J$-holomorphic maps $ \mul{u}_S$ produced by Theorem
\ref{gluingstable} with parameter space a neighborhood $S$ of $0$ in
$\Def(\mul{u}),$ equipped with a canonical identification $\iota$ of
the central fiber with the original map $\mul{u}$.  In the case that
the domain $\mul{\Sigma}$ is not a stable (marked) curve, we choose
codimension two submanifolds $\mul{Y} = (Y_1,\ldots, Y_l)$ meeting
$\mul{u}$ transversally so that $\mul{\Sigma}$ with the additional
marked points is stable.  Applying this to the family $\mul{u}_S$ gives
a family of marked stable maps $\mul{u}_S^{\mul{Y}}$ with $n + l$ marked
points over a parameter space $S \subset \Def(\mul{u}^{\mul{Y}})$ in the
deformation space of the map with the additional marked points.  Now
$\Def(\mul{u}^{\mul{Y}}) \cong \Def(\mul{u}) \oplus \bigoplus_{i=1}^l
T_{z_i} \mul{\Sigma}$ includes the deformations of the markings, but
these are fixed by requiring that the additional marked points map to
the given collection $\mul{Y}$.  Forgetting the additional marked
points gives a family $\mul{u}_S$ of stable maps with $n$ marked points
over a neighborhood of $0$ in $\Def(\mul{u})$.

\begin{proposition}  $(\mul{u}_S,\iota)$ is a versal stratified-smooth 
deformation of $\mul{u}$, and in fact $\mul{u}_S$ gives a versal
stratified-smooth deformation of any of its fibers.
\end{proposition} 

\begin{proof}  First suppose that $\mul{\Sigma}$ is stable.   
Let $(\mul{u}_{S^1}^1,\iota^1)$ be another stratified-smooth
deformation of $\mul{u}$ with parameter space $S^1$.  Let
$\mul{\Sigma}_S \to S \subset \Def(\mul{\Sigma})$ be a minimal versal
deformation of $\mul{\Sigma}$.  The family $\mul{\Sigma}_{S^1}^1$ is
obtained by pull-back of $\mul{\Sigma}_S$ by a stratified-smooth map
$\psi: S^1 \to S$.  By definition the map $\mul{u}_s^1$ converges to
the central fiber in the Gromov topology as $s$ converges to the base
point $0 \in S^1$.  The exponential decay estimate of \cite[Lemma
  4.7.3]{ms:jh} for holomorphic cylinders of small energy imply that
for $s$ sufficiently close to $0$, $\mul{\Sigma}_s^1,\mul{u}_s^1$ is
given by exponentiation, $ \mul{u}_s^1 =
\exp_{\mul{u}^{\mul{\delta}}}(\mul{\xi}) $ for some $\mul{\xi} \in
\Omega^0( \mul{u}^{\mul{\delta},*} TX)$ with $\Vert \mul{\xi}
\Vert_{1,p} < \eps_1$, for $s$ sufficiently close to $0$.  Proposition
\ref{surject} produces a stratified-smooth map $\psi : S^1 \to
\Def(\mul{u})$ such that $\mul{u}_{S^1}^1$ is the pull-back of $\psi$.
To show that the deformation $(\mul{u}_S, \iota)$ is universal, let $
\phi_j: \mul{\Sigma}_{S^1}^1 \to \psi_j^* \mul{\Sigma}_S, j = 0,1 $ be
isomorphisms of families inducing the identity on the central fiber.
The difference between the two automorphisms is an automorphism of the
family $\mul{\Sigma}_{S^1}^1$ inducing the identity on the central
fiber; since the automorphism group of the central fiber is discrete,
the automorphism must be the identity.  In the case that
$\mul{\Sigma}$ is not stable, after adding marked points passing
through $Y_1,\ldots, Y_l$, we obtain a family
$\mul{u}_{S^1}^{1,\mul{Y}}$ of stable maps with $n+l$ marked points.
By the case with stable domain, this family is obtained by pull-back
of $\mul{u}_S^{\mul{Y}}$ by some map $S^1 \to S$.  Hence
$\mul{u}_{S^1}^1$ is obtained by pull-back by the same map.  The
argument for an arbitrary fiber is similar and left to the reader.
\end{proof} 

\begin{remark}  In the case that $\mul{\Sigma}$ is unstable, it seems
likely that restricting the family of Theorem \ref{gluingstable} to
$\Def(\mul{u})$ (that is, the perpendicular of $\aut(\mul{\Sigma})$)
also gives a universal deformation, but we do not know how to prove
this.  The problem is that in this case, several different gluing
parameters give the same curve, and we do not have an implicit
function theorem for varying gluing parameter.
\end{remark} 

\noindent {\em Step 7: Injectivity}

By injectivity, we mean that the family constructed above contains
each nearby stable map exactly once, up to the action of $\Aut(u)$.
This is part of what we called ``strongly universal'' in Definition
\ref{excellent}.

\begin{theorem}  The versal deformations constructed in Step 6 above 
are strongly universal.
\end{theorem}

\begin{proof}  Let $\mul{u}_S$ be a deformation constructed
as in Step 6, using the exponential gluing profile.  Let
$\mul{\Sigma}_{1,S^1} \to S^1$ be a family giving a universal
deformation of the curve $\mul{\Sigma}^{\mul{Y},\mul{u}}$ obtained by
adding the additional markings mapping to the given submanifolds.  By
Definition \ref{compat}, the family $\mul{\Sigma}_{S}^{\mul{Y},\mul{u}}$
induces a map $\phi: S \to S_1$ whose differential is injective in a
neighborhood of $0$.  By the inverse function theorem for $C^1$ maps,
$\phi$ induces a homeomorphism onto its image.  In particular, any two
distinct fibers of $\mul{\Sigma}_S^{\mul{Y},\mul{u}}$ are non-isomorphic,
and so two fibers of $\mul{\Sigma}_S$ are isomorphic if and only if
they are related by a permutation of the markings.  After shrinking
$S$, this happens only if the permutation is induced by an
automorphism of $\mul{u}$.  Given another family $\mul{u}_{S'}' :
\mul{\Sigma}_{S'} \to S'$ corresponding to a deformation of a fiber of
$\mul{u}_S \to S$, by the uniqueness part of the implicit function
theorem, a map $\phi': S' \to S_1$ so that $\mul{u}_{S'}'$ is obtained
by pull-back from $\mul{u}_S$, and this map is unique by the
injectivity just proved.  This shows that $\mul{u}_S$ gives a
stratified-smooth universal deformation of any of its fibers, and so
is strongly universal.
\end{proof}  

The Theorem implies that the families in the universal deformations
constructed above define stratified-smooth-compatible charts for the
moduli space $\ol{M}_{g,n}(X,d)$.  That is, for any stratum
$M_{g,n,\Gamma}(X,d)$, the restriction of the charts given by the
universal deformation of some map of type $\Gamma$ to
$M_{g,n,\Gamma}(X,d)$ are smoothly compatible.

\begin{corollary} \label{orbi}
Let $X,J$ be as above.  For any $g \ge 0, n \ge 0$, the strongly
universal stratified-smooth deformations of parametrized regular
stable maps provide $\ol{M}_{g,n}^{\reg}(X)$ with the structure of a
stratified-smooth topological orbifold.
\end{corollary} 

In order to apply localization one needs to know that the fixed point
sets admit tubular neighborhoods. For this it is helpful to know that
$\ol{M}_{g,n}^\reg(X,d)$ admits a $C^1$ structure.  In order to obtain
compatible charts, we construct the local coordinates inductively as
in Definition \ref{compatcoord}, starting with the strata of highest
codimension.

\begin{proposition}  \label{diffable} 
Let $X,J$ be as above.  For any compatible system of local coordinates
near the nodes, the strongly universal deformations constructed using
the exponential gluing profile equip $\ol{M}_{g,n}^\reg(X)$ with the
structure of a $C^1$-orbifold.
\end{proposition} 

\begin{proof}  We claim that the charts induced by the universal
deformations are $C^1$-compatible, assuming they are constructed from
the same system of local coordinates near the nodes.  Given two sets
of submanifolds $\mul{Y}_1,\mul{Y}_2$, define $\mul{Y} = \mul{Y}_1 \cup
\mul{Y}_2$.  The family $ \mul{\Sigma}^{\mul{Y},\mul{u}}$ admits proper
\'etale forgetful maps $ \mul{\Sigma}^{\mul{Y},\mul{u}}_S \to
\mul{\Sigma}_S^{\mul{Y}_j,\mul{u}} , \quad j= 1,2 .$ The fiber consists
of reorderings of the additional marked points induced by the action
of $\Aut(\mul{\Sigma},\mul{u})$, and the diagram provided by
$\mul{\Sigma}^{\mul{Y},\mul{u}}$ expresses the composition as a smooth
$C^1$-morphism of orbifolds.
\end{proof}  

\begin{remark}  \label{smooth}  Any compact $C^1$ orbifold admits a compatible 
$C^\infty$ structure, in analogy with the situation with manifolds.
  Indeed, as is well known any orbifold admits a presentation as the
  quotient of a manifold (namely its orthogonal frame bundle) by a
  locally free group action, and so the orbifold case follows from the
  equivariant case proved in Palais \cite{pal:act}.  Hence
  $\ol{M}^{\reg}_{g,n}(X,d)$ if compact admits a (non-canonical)
  smooth structure.  Presumably the compactness assumption may be
  removed but we have not proved that this is so.  See however the
  construction of smoothly compatible Kuranishi charts in
  \cite[Appendix]{fooo}.
\end{remark} 

\section{Deformations of symplectic vortices}

We begin by reviewing the theory of symplectic vortices introduced by
Mundet i Riera \cite{mun:ham} and Salamon and collaborators
\cite{ciel:vor}.  Let $\Sigma$ be a compact complex curve, $G$ a
compact Lie group, and $\pi: P \to \Sigma$ a smooth principal
$G$-bundle.  Given any left $G$-manifold $F$ we have a left action of
$G$ on $P \times F$ given by $g(p,f) = (pg^{-1},gf)$ and we denote by
$P(F) = (P \times F)/G$ the quotient, that is, the associated fiber
bundle with fiber $F$.  Let $X$ be a compact Hamiltonian $G$-manifold
with symplectic form $\omega$ and moment map $\Phi:X \to \g^*$.  The
action of $G$ on $X$ induces an action on $\J(X)$; and we denote by
$\J(X)^G$ the invariant subspace. Let $\psi: \Sigma \to BG$ be a
classifying map for $P$, so that $P \cong \psi^* EG$ and $P(X) \cong
\psi^* EG \times_G X \cong \psi^* X_G$ where $X_G = EG \times_G X$.
Continuous sections $u: \Sigma \to P(X)$ are in one-to-one
correspondence with lifts of $\psi$ to $X_G$.  The homology class
$\deg(u)$ of the section $u$ is defined to be the homology class
$\deg(u) \in H_2^G(X,\Z)$ of the corresponding lift.  Let $\A(P)$ be
the space of smooth connections on $P$, and $P(\g)$ the adjoint
bundle.  For any $A \in \A(P)$, let $F_A \in \Omega^2(\Sigma,P(\g))$
the curvature of $A$.  Any connection $A \in \A(P)$ induces a map of
spaces of almost complex structures
$$ \J(X)^G \to \J(P(X)), \ \ J \mapsto J_A$$
by combining the almost complex structures on $X$ and $\Sigma$ using
the splitting defined by the connection.  Let $\Gamma(\Sigma,P(X))$
denote the space of smooth sections of $P(X)$.  Consider the vector
bundle
\begin{equation} \label{crop}
\bigcup_{u \in \Gamma(\Sigma,P(X)) } \Omega^{0,1}(\Sigma,
u^*T^{\on{vert}} P(X)) \to \Gamma(\Sigma,P(X)) .
\end{equation}
We denote by $\olp_A$ the section given by the Cauchy-Riemann operator
defined by $J_A$.  A {\em gauged map} from $\Sigma$ to $X$ is a datum
$(P,A,u)$ where $A \in \A(P)$ and $u: \Sigma \to P(X)$ is a section.
A {\em gauged holomorphic map} is a gauged map $(P,A,u)$ such that
$\olp_A u = 0 $.  Let $\H(P,X)$ be the space of gauged holomorphic
maps with underlying bundle $P$.  Let $\G(P)$ denote the group of
gauge transformations
$$ \G(P) = \{ a: P \to P,  a(pg) = a(p)g, \ \ \ \pi \circ a = \pi \} .$$
The Lie algebra of $\G(P)$ is the space of sections
$\Omega^0(\Sigma,P(\g))$ of the adjoint bundle $P(\g) = P \times_G
\g$.  We identify $\g \to \g^*$, and hence $P(\g) \to P(\g^*)$, using
an invariant metric on $\g$.  Let $P(\Phi): P(X) \to P(\g)$ denote the
map induced by the equivariant map $\Phi : X \to \g$.

\begin{definition} 
A gauged holomorphic map $(A,u) \in \H(P,X)$ is a {\em symplectic
  vortex} (or vortex for short) if it satisfies
$$ F_A + \Vol_{\Sigma} u^* P(\Phi) = 0 .$$
An {\em $n$-marked} symplectic vortex is a vortex $(A,u)$ together
with $n$-tuple $\mul{z} = (z_1,\ldots, z_n)$ of distinct points on
$\Sigma$.  A marked vortex $(A,u,\mul{z})$ is {\em stable} if it has
finite automorphism group.
\end{definition} 

The equation in the definition can be interpreted as the zero level
set condition for a formal moment map for the action of the group of
gauge transformations, see \cite{mun:ham}, \cite{ciel:vor}.  The {\em
  energy} of a gauged holomorphic map $(A,u)$ is given by
$$ E(A,u) = \hh \int_\Sigma \left(| \d_A u |^2 + |F_A|^2 + |u^*
 P(\Phi)|^2 \right)  \Vol_\Sigma.$$
The \emph{equivariant symplectic area} of a pair $(A,u)$ is the
pairing of the homology class $\deg(u)$ with the class $[\omega_G =
  \omega + \Phi] \in H^2(X_G)$,
$$ D(A,u) = (\deg(u), [\omega_G]) = ([\Sigma], u^* [\omega_G] ) .$$

\begin{lemma}  Suppose $\Vol_\Sigma$ is the K\"ahler form 
for the metric on $\Sigma$.  The energy and equivariant area are
related by
\begin{equation}  \label{energyaction} 
E(A,u) = D(A,u) + \int_\Sigma \left( | \olp_A u |^2 + \hh | F_A + \Vol_\Sigma
u^* P(\Phi) |^2 \right) \Vol_\Sigma.
\end{equation}
\end{lemma} 

\begin{proof}  See \cite[Proposition 2.2]{ci:symvortex}. 
\end{proof} 

In particular, for any symplectic vortex the energy and action are
equal.  Let $M(P,X)$ denote the moduli space of vortices
$$ M(P,X) := \H(P,X) \qu \G(P) = \{ F_A + \Vol_{\Sigma} u^* P(\Phi) =
0 \}/\G(P) .$$
Let $M_n(P,X)$ denote the moduli space of $n$-marked vortices, up to
gauge transformation, and
$M_n(\Sigma,X) = \bigcup_{P \to \Sigma} M_n(P,X) $
the union over types of bundles $P$.  Clearly, $M_n(\Sigma,X)$ is
homeomorphic to the product $M(\Sigma,X) \times M_n(\Sigma) $ where
$M(\Sigma,X) := M_0(\Sigma)$ and $M_n(\Sigma)$ denotes the
configuration space of $n$-tuples of distinct points on $\Sigma$.

We wish to study families and deformations of symplectic vortices.
For families with smooth domain, the definitions are straightforward:

\begin{definition}  \label{strongly universal}
A {\em smooth family of vortices} on a principal $G$-bundle $P$ on
$\Sigma$ over a parameter space $S$ consists of a family of
connections depending smoothly on $s \in S$, that is, a smooth map
$A_S : S \times P \to T^*P \otimes \g$ on $P$ such that the
restriction $A_s$ of $A_S$ to any $ \{s \} \times P$ is a connection,
together with a smooth family of (pseudo)holomorphic sections $u_S =
(u_s)_{s \in S}$, such that each pair $(A_s,u_s), s \in S$ is a
symplectic vortex.  A {\em deformation} of $(A,u)$ is a germ of a
smooth family $(A_S,u_S)$ together with an isomorphism (gauge
transformation) relating $(A_0,u_0)$ with $(A,u)$.  A deformation is
{\em universal} if it satisfies the condition as in Definition
\ref{ssfam}, and {\em strongly universal} if it satisfies the
conditions in Definition \ref{excellent}. \end{definition}

We define a {\em linearized operator} associated to a vortex as
follows.  Define
\begin{equation} \label{linearized1}
 \d_{A,u}: \Omega^1(\Sigma,P(\g)) \oplus \Omega^0(\Sigma,u^*
 TP(X) ) \to \Omega^2(\Sigma,P(\g)) \end{equation}
$$ \d_{A,u}(a,\xi) := \d_A a + \Vol_{\Sigma} u^* L_\xi
P(\Phi).
$$
Here $L_\xi P(\Phi)$ denotes the derivative of $P(\Phi)$ with respect
to the vector field generated by $\xi$, and $u^* L_\xi P(\Phi)$ its
evaluation at $u$.  Define an operator
\begin{equation} \label{linearized2} 
 \d_{A,u}^*: \Omega^1(\Sigma,P(\g)) \oplus \Omega^0(\Sigma,u^* TP(X) )
 \to \Omega^0(\Sigma,P(\g)) \end{equation}
$$ \d_{A,u}^*(a,\xi) = \d_A^*a + u^* L_{J \xi} P(\Phi) .$$
(This is not the adjoint of operator in \eqref{linearized}, but rather
defined by analogy with the case $X$ trivial.)  It is shown in
\cite[Section 4]{ciel:vor} that if $(A,u)$ is stable then the set
$$ W_{A,u} = \{  (A + a, \exp_u(\xi)), (a,\xi) \in \ker \d_{A,u}^*\} $$
is a slice for the gauge group action near $(A,u)$.  Define
\begin{multline} \label{Feps}
 \cF_{A,u} : \Omega^1(\Sigma,P(\g)) \oplus
\Omega^0(\Sigma,u^* T^{\on{vert}} P(X)) \\ \to (\Omega^0 \oplus
\Omega^2)(\Sigma,P(\g)) \oplus \Omega^{0,1}(\Sigma,u^*
T^{\on{vert}} P(X)) \\
 (a,\xi) \mapsto \left( F_{A + a} +  \Vol_\Sigma
\exp_u(\xi)^* P(\Phi), \d_{A,u}^* (a,\xi), \Psi_u(\xi)^{-1}
\ol{\partial}_{A + a} \exp_u(\xi) \right) .\end{multline}
Let 
$$\Omega^1(\Sigma,P(\g)) \to \Omega^1(\Sigma,u^*
T^{\on{vert}}P(X)), \ \ \ a \mapsto a_X $$
denote the map induced by the infinitesimal action.  The linearization
of the last component \eqref{Feps} is
$$ D_{A,u} (a,\xi) = (\nabla_A \xi)^{0,1} + \hh J_u (\nabla_\xi J )_u \partial_A u
+ a_X^{0,1}.$$
Here $0,1$ denotes projection on the $0,1$-component.  The linearized
operator for a vortex $(A,u)$ is the operator
\begin{multline}  \label{linearized}
 \ti{D}_{A,u} = (\d_{A,u}, \d_{A,u}^*, D_{A,u} ):
\Omega^1(\Sigma,P(\g)) \oplus \Omega^0(\Sigma, u^* T^{\on{vert}}P(X))
 \\ \to (\Omega^0 \oplus \Omega^2)(\Sigma,P(\g)) \oplus 
\Omega^{0,1}(\Sigma,  u^* T^{\on{vert}}P(X))
\end{multline}
A vortex $(A,u)$ is {\em regular} if the operator $\ti{D}_{A,u}$ is
surjective.  A marked vortex $(A,u,\mul{z})$ is regular if the
underlying unmarked vortex is regular.  The {\em space of
infinitesimal deformations} of $(A,u)$ is
$ \Def(A,u) := \ker(\ti{D}_{A,u}) .$

\begin{theorem}  \label{deformvortex} 
Any regular vortex with smooth domain $(A,u)$ has a strongly universal
smooth deformation if and only if it is stable.
\end{theorem}  

\begin{proof} 
Give the spaces of connections and sections the structure of Banach
manifolds by taking completions with respect to Sobolev norms $1,p$
for $1$-forms, and $0,p$ for $0$ and $2$-forms. For $p>2$, the map
$\cF_ {A,u}$ is a smooth map of Banach spaces. 
\begin{multline} 
 \cF_{A,u} : \Omega^1(\Sigma,P(\g))_{1,p} \oplus \Omega^0(\Sigma,u^*
 T^{\on{vert}} P(X))_{1,p} \\ \to (\Omega^0 \oplus
 \Omega^2)(\Sigma,P(\g))_{0,p} \oplus \Omega^{0,1}(\Sigma,u^*
 T^{\on{vert}} P(X))_{0,p} \end{multline}
equivariant for the action of the group $\G(P)_{2,p}$ of gauge
transformations of class $2,p$.  Suppose that $(A,u)$ is regular and
stable.  By the implicit function theorem, there is a local
homeomorphism
$$ \ker(\ti{D}_{A,u}) \to \left\{  \begin{array}{c} F_{A
    + a} + \Vol_\Sigma (\exp_u(\xi))^* P(\Phi) = 0 \\ \olp_{A + a}
  (\exp_u(\xi)) = 0 \\ \ \d_{A,u}^*(a,\xi) = 0 \end{array} \right\}
.$$
This gives rise to a family $(A_S,u_S) \to S$ over a neighborhood $S$
of $0$ in $\ker(\ti{D}_{A,u})$.  By \cite[Theorem 3.1]{ciel:vor},
$(A_S,u_S)$ is a smooth family, assuming $(A,u)$ is smooth.  Given any
other family $(A_{S'}',u_{S'}') \to S'$ of stable vortices with
$(A'_0,u_0') = (A,u)$, the implicit function theorem provides a smooth
map $S' \to S$ so that $(A'_{S'},u'_{S'})$ is obtained from $(A,u)$ by
pull-back.  The first property of the universal deformation is a
consequence of the slice condition; the second property follows from
the fact that the projection $\ker(\ti{D}_{A,u}) \to
\ker(\ti{D}_{A_s,u_s})$ is an isomorphism for sufficiently small $s$.
\end{proof}  
\noindent Let $M^\reg_n(\Sigma,X)$ denote the moduli space of regular,
stable $n$-marked symplectic vortices from $\Sigma $ to $X$.  We
denote by $(c_1^G(TX),d)$ the pairing of $d$ with the first Chern
class $c_1^G( P(TX) \to P(X))$

\begin{theorem}  \label{smoothreg}  Let $\Sigma,X,J$ be as above.  
 $M^\reg_n(\Sigma,X)$ has the structure of a smooth orbifold with
  tangent space at $[A,u]$ isomorphic to $\Def(A,u)$, and dimension of
  the component of homology class $d \in H_2^G(X)$ is given by
$$
 \dim(M^\reg_n(\Sigma,X,d))
= (1-g)(\dim(X) - 2\dim(G)) +
2((c_1^G(TX),d) + n) .$$
\end{theorem}  

\begin{proof} 
Charts for $M^{\reg}_n(\Sigma,X)$ are provided by the strongly
universal deformations.  The dimension of the tangent space at $[A,u]$
is given by the index of the linearized operator $\ti{D}_{A,u}$, which
deforms via Fredholm operators to the sum of the operator $\d_A \oplus
\d_A^*$ for the connection, which has index $2\dim(G)(g-1)$, and the
linearized Cauchy-Riemann operator on the nodal curve, which has index
$(1-g)\dim(X) + 2n + 2(c_1^G(TX),d)$ by Riemann-Roch, if $(A,u)$ has
equivariant homology class $d$ (which determines the first Chern class
of $P$ by projection.)
\end{proof}  

\subsection{Polystable vortices}

The moduli space of symplectic vortices admits a compactification
which allows bubbling of the section in the fibers. 

\begin{definition} \label{nodal} 
A {\em nodal gauged marked holomorphic map} from $\Sigma$ to $X$
consists of a datum $(\hat{\Sigma},P,A,\mul{u},\mul{z})$ where $P \to
\Sigma$ is a principal $G$-bundle, $A \in A(P)$ is a connection,
$\hat{\Sigma}$ is a marked nodal curve, $v: \hat{\Sigma} \to \Sigma$
is a holomorphic map of degree $[\Sigma]$, and $\mul{u}: \hat{\Sigma}
\to P(X)$ is a $J_A$-holomorphic map from a nodal curve $\hat{\Sigma}$
such that $\pi \circ \mul{u}$ has class $[\Sigma]$.  In other words,
\begin{enumerate}
\item $\hat{\Sigma}$ is a connected nodal complex curve consisting of
  a {\em principal component} $\Sigma_0$ equipped with an isomorphism
  with $\Sigma$ together with a number of projective lines
  $\Sigma_1,\ldots,\Sigma_k$.  We denote by $w_1^\pm,\ldots,w_m^\pm$
  the nodes.  For each $i = 1,\ldots, k$, we denote by $w_i^0 \in
  \Sigma_0$ the attaching point to the principal component.
\item 
$(A,u) \in \H(P,X)$ 
is a gauged holomorphic map from $\Sigma$ to $X$;
\item for each non-principal component $\Sigma_i$, a holomorphic map
  $u_i: \Sigma_i \to P(X)_{w_i^0}$;
\item $\mul{z} = (z_1,\ldots,z_n) \in \hat{\Sigma}$ are distinct, smooth
points of $\hat{\Sigma}$.
\end{enumerate}
\end{definition} 

Let $\H(\hat{\Sigma},P(X))$ denote the space of nodal gauged
holomorphic sections with domain $\hat{\Sigma}$ and bundle $P$. The
group of gauge transformations $\G(P)$ acts on $\H(\hat{\Sigma},P(X))$
by $ g (A,\mul{u}) = (g^* A, g \circ \mul{u}) .$ The generating vector
field for $\zeta \in \Omega^0(\Sigma,P(\g))$ acting on
$\H(\hat{\Sigma},P(X))$ at $(\hat{\Sigma},A,\mul{u})$ is the tuple
given by
\begin{equation} \label{generate}
\zeta_{\H(\hat{\Sigma},P(X))}(\hat{\Sigma},A,\mul{u}) = (\d_A \zeta,
u_0^* P(\zeta_X), (u_i^* P(\zeta_X(w_i^0)))_{i=1}^k )
 \end{equation}
in $ \Omega^1(\Sigma,P(\g)) \oplus \Omega^0(\hat{\Sigma},\mul{u}^*
T^{\on{vert}} P(X))$.  Here $P(\zeta_X) \in \Omega^0(\Sigma, P(
\Vect(X)))$ is the fiber-wise vector field generated by $\zeta$ and
$u_0^* P(\zeta) \in T^{\on{vert}} P(X)$ is the evaluation at $u_0$.
Similarly for the bubble components $u_1,\ldots, u_k$ in the fibers
$P(X)_{w^0_1},\ldots, P(X)_{w^0_k}$.  A slice is given by taking the
perpendicular to the tangent spaces to the $\G(P)$-orbits.  We will
assume for simplicity that the stabilizer of the $\G(P)$ action on the
principal component is finite, so that a slice is given locally by the
kernel of $\d_{A,u_0}^*$, that is, the Coulomb gauge condition on the
principal component.  The implicit function theorem shows that any
nearby pair $(A_1,\mul{u}_1)$ is complex gauge equivalent to a pair of
the form $(A + a, \exp_{\mul{u}}(\mul{v}))$ with $(a,\mul{v}) \in \ker
\d_{A,u_0}^*$.

\begin{definition}
A {\em nodal vortex} is a stable nodal gauged holomorphic map such
that the principal component is an vortex.  A nodal vortex
$(\hat{\Sigma},A,\mul{u},\mul{z})$ is {\em polystable} if each sphere
bubble $\Sigma_i$ on which $u_i$ is constant has at least three marked
or singular points, and {\em stable} if it has finite automorphism
group.
 An {\em isomorphism} of nodal vortices
 $(\hat{\Sigma},A,\mul{u},\mul{z}),(\hat{\Sigma}',A',\mul{u}',\mul{z}')$
 consists of an automorphism of the domain, acting trivially on the
 principal component, and a corresponding automorphism of the
 principal bundle mapping $(A,\mul{u}) $ to $(A',\mul{u}')$ and mapping
 the markings $\mul{z}$ to $\mul{z}'$.
For any nodal section $\mul{u}:\hat{\Sigma} \to P(X)$, the {\em homology
  class} of $\mul{u}$ is defined as the sum of the homology class $d_0
\in H_2^G(X,\Z)$ of the principal component $u_0$ and the homology
classes $d_i \in H_2(X,\Z), i = 1,\ldots,k$ of the sphere bubbles,
using the inclusion $H_2(X,\Z) \to H_2^G(X,\Z)$ given by equivariant
formality.
The {\em combinatorial type} $\Gamma(\hat{\Sigma},A,\mul{u},\mul{z})$ of
a gauged nodal map is a rooted graph whose vertices represent the
components of $\hat{\Sigma}$, whose finite edges represent the nodes,
semi-infinite edges represent the markings, and whose root vertex
represents the principal component.
\end{definition} 

\noindent 
Note that there is no condition for points on the principal
component. In particular, nodal gauged holomorphic maps with no
markings can be polystable.  The term {\em polystable} is borrowed
from the vector bundle case. In that situation, a bundle is {\em
  stable} if it is flat and has only central automorphisms and
polystable if it is a direct sum of stable bundles of the same slope.
Any flat bundle is automatically polystable; a bundle is {\em
  semistable} if it is grade equivalent to a polystable bundle. In
particular, the moduli space of stable bundles is definitely {\em not}
compact, and we feel that the vortex terminology should include this
fact as a special case.

\noindent From now on, we fix the bundle $P$.

\begin{definition}  Let $X$ as above.  
A {\em smooth family of fixed type} of nodal vortices to $X$ consists
of a smooth family $\hat{Sigma}_S \to S$ of nodal curves of fixed
type, a smooth family of holomorphic maps $v_S: \hat{\Sigma}_S \to
\Sigma$ of class $[\Sigma]$, a smooth family $\mul{u}_S :
\hat{\Sigma}_S \to P(X)$ of maps, and a smooth family $A_S : S \times
P \to T^* P$ of connections over $S$.  A {\em smooth deformation} of a
nodal vortex $(A,\hat{\Sigma},\mul{u},\mul{z})$ of fixed type consists
of a germ of a smooth family
$(A_S,\hat{\Sigma}_S,\mul{u}_S,\mul{z}_S)$ of nodal vortices of fixed
type together with an identification $\iota$ of of the central fiber
with $(A,\hat{\Sigma},\mul{u},\mul{z})$.  A {\em stratified-smooth
  family of marked nodal symplectic vortices} is a datum
$(\hat{\Sigma}_S,A_S,\mul{u}_S,\mul{z}_S)$ consisting of a
stratified-smooth family $\hat{\Sigma}_S \to S$ of nodal curves, a
stratified-smooth family of holomorphic maps $v: \hat{\Sigma}_S \to
\Sigma$ of class $[\Sigma]$, a stratified-smooth family $A_S$ of
connections on $P$, a stratified-smooth family of maps $\mul{u}_S:
\hat{\Sigma}_S \to P(X)$; such that each triple
$(\hat{\Sigma}_s,A_s,\mul{u}_s,\mul{z}_s)$ is a marked nodal
symplectic vortex.  A {\em family of polystable symplectic vortices}
is a family of marked nodal symplectic vortices such that any fiber is
polystable.
\end{definition}  

\noindent A smooth vector bundle $\mul{E} \to \hat{\Sigma}$ is a
collection of smooth vector bundles $E_i$ over the components
$\Sigma_i$ of $\hat{\Sigma}$, equipped with identifications of the
fibers at nodal points $E_{i^+(w^+_k)} \to E_{i^-(w^-_k)}, k =
1,\ldots, m$.  We denote by $\Omega(\hat{\Sigma},\mul{E})$ the sum
over components, $ \Omega(\hat{\Sigma},\mul{E}) = \bigoplus_{i=1}^k
\Omega(\Sigma_i, E_i) $ where $E_i = \mul{E} | \Sigma_i .$

\begin{definition}  For a polystable vortex
$(\hat{\Sigma},A,\mul{u})$, let $\ti{D}_{A,\mul{u}}$ denote the {\em
    linearized operator}
\begin{multline}
 \Omega^1({\Sigma},P(\g)) \oplus 
\Omega^0(\hat{\Sigma},\mul{u}^*
 T^{\on{vert}} P(X))
 \\ \to (\Omega^{0}
 \oplus \Omega^2)({\Sigma},P(\g)) \oplus \Omega^{0,1}(\hat{\Sigma},\mul{u}^*
 T^{\on{vert}} P(X)) \oplus \bigoplus_{k=1}^m T_{\mul{u}(w_k^\pm)}^{\on{vert}} P(X)
\end{multline}
given by the linearized vortex operator $(\d_{A,u_0}, D_{A,u_0})$ on
the principal component, the linearized Cauchy-Riemann operator
$\ti{D}_{u_i}$ on the bubbles, the slice operator $\d_{A,u_0}^*$, and
the difference operator on the fibers over the nodes
$$ \Omega^0(\hat{\Sigma},\mul{u}^* T^{\on{vert}} P(X)) \to
\bigoplus_{i=1}^m T_{\mul{u}(w_i^\pm)}^{\on{vert}} P(X), \quad \mul{\xi}
\mapsto ( \mul{\xi}(w_i^+) - \mul{\xi}(w_i^-))_{i=1}^m .$$
$(A,\mul{u})$ is {\em regular} if $\ti{D}_{A,\mul{u}}$ is surjective.
The space of infinitesimal deformations of $A,\mul{u}$ of fixed type is
$ \Def_\Gamma(A,\mul{u}) := \ker \ti{D}_{A,\mul{u}}/\aut(\hat{\Sigma})
$
where $\aut(\hat{\Sigma})$ denotes the group of infinitesimal
automorphisms acting trivially on the principal component.  The space
of infinitesimal deformations of $A,\mul{u}$ is
$$ \Def(A,\mul{u}) := \Def_\Gamma(A,\mul{u}) \oplus 
\bigoplus_{i=1}^m
 T_{w_j^+} \hat{\Sigma} \otimes T_{w_j^-} \hat{\Sigma} $$
consisting of a deformation of fixed type together with a collection
of gluing parameters.
\end{definition} 

\subsection{Constructing deformations of symplectic vortices}

First we construct deformations of fixed type.

\begin{theorem}  \label{deformnodalvortex} 
A regular polystable vortex has a strongly universal smooth
deformation of fixed type if and only if it is stable.
\end{theorem}  

\noindent The proof is by the implicit function theorem applied to the
map
\begin{multline} \cF_{A,\mul{u}} (a,\mul{\xi}) = (F_{A + a} + \Vol_\Sigma
(\exp_{u_0}(\xi_0))^* P(\Phi), \d_{A,u_0}^* (a,\mul{\xi}), \\
  \Psi_{u_0}(\xi_0)^{-1} \olp_A u_0, (\Psi_{u_i}(\xi_i)^{-1} \olp
  u_i)_{i=1}^k, (\mul{\xi}(w_i^+) - \mul{\xi}(w_i^-))_{i=1}^m )
\end{multline}
whose linearization is $\ti{D}_{A,\mul{u}}$.  The proof is left to the
reader.  We denote by $M_{n,\Gamma}(\Sigma,X,d)$ of the moduli space
of isomorphism classes of polystable vortices of combinatorial type
$\Gamma$ of homology class $d \in H_2^G(X,\Z)$, and
$M^{\reg}_{n,\Gamma}(\Sigma,X,d)$ the regular locus.

\begin{corollary}  \label{stablemreg}  
$M_{n,\Gamma}^\reg(\Sigma,X,d)$ has the structure of a smooth orbifold
  of dimension given, if $\Sigma$ is connected, by $(1 - g)(\dim(X) -
  2\dim(G)) + 2( n + (c_1^G(TX),d) - m)$ where $m$ is the number of
  nodes.
\end{corollary} 

\label{gluingv} 

We now prove that a regular stable symplectic vortex from $\Sigma$ to
$X$ admits a strongly universal stratified-smooth deformation if it is
strongly stable, that is, Theorem \ref{main}.  We explain the
construction for a single bubble only, so that $\hat{\Sigma}$ is the
union of a principal component $\Sigma_+ = \Sigma$ and a holomorphic
sphere $\Sigma_-$, attached by a single pair $w_\pm$ of nodes.  We
denote by $(A,u_+)$ the restriction to the principal component and by
$u_-$ the bubble, so that $x := u_+(w_+) = u_-(w_-)$ and $\mul{u} =
(u_+,u_-)$.  We choose a local coordinate near $w$, equivariant for
the action of the automorphism group $\Aut(A,\mul{u})$ in the sense
that $\Aut(A,\mul{u})$ acts on the local coordinate by multiplication
by roots of unity.  The construction depends on the following choices:

\begin{definition}  A {\em gluing datum} for $(\hat{\Sigma},A,\mul{u})$ consists of 
\begin{enumerate} 
\item neighborhoods $U_\pm$ of the nodes $w_\pm$;
\item {\em local coordinates} $\mul{\kappa} = (\kappa_+,\kappa_-)$ on
  $U_\pm$;
\item a trivialization $ P |_{U_0} \to G \times U_0 $ on $U_0$;
\item a {\em gluing parameter} $\delta$;
\item an {\em annulus parameter} $\rho$ ;
\item a {\em cutoff function} $\alpha$ as in \eqref{firstcut}. 
\end{enumerate}
\end{definition} 

\noindent {\em Step 1:  Approximate Solution} 

 Given a nodal vortex $(A,u)$ as above and a gluing datum we wish to
 define an approximate solution to the vortex equations $(A,
 u^\delta)$.  Let $\exp_x: T_x X \to X $ denote the exponential map
 defined by the metric on $X$.  Define sections
$$\xi_\pm: U_\pm \to T_x X, \quad u(z)=\exp_{x}(\xi_\pm(z)) .$$
Let $\hat{\Sigma}^\delta$ denote the surface obtained by gluing; since
the bubble is genus zero, this surface is isomorphic to $\Sigma$ but
not canonically.  Define the {\em pre-glued section} $\mul{u}^\delta:
\hat{\Sigma}^\delta \to P(X)$,
\begin{multline} 
\mul{u}^\delta(z) = \exp_{x} \left( \alpha(|\kappa_+(z)|/ \rho
|\delta|^{1/2}) (\xi_+(z) - \xi_\pm(w^\pm)) \right. +\\  \left. \alpha
(|\kappa_-(z)|/ \rho |\delta|^{1/2} ) (\xi_-(z) - \xi_\pm(w^\pm)) +
\xi_\pm(w^\pm) \right)
\end{multline} 
for $| \kappa_\pm(z) | \leq 2|\delta|^{1/2} \rho^2$; elsewhere let
$\mul{u}^\delta(z) = \mul{u}(z)$, using the identification of $\Sigma$
with $\hat{\Sigma}$ away from the gluing region.  We do not modify $A$
in the bubble region; this is because after re-scaling the connection
on the bubble is already close to the trivial connection.

The pair $(A,\mul{u}^\delta)$ is an approximate solution to the vortex
equations in the following Sobolev norms.  Let
$\Omega(\Sigma^\delta,P^\delta(\g))_{k,p}$ denote the
$W^{k,p}$-completion using the standard metric. Let $g^\delta$ denote
the $C^0$ metric on the glued surface in \eqref{gluemetric}.  Let
$\Omega^{0,1}(\Sigma^\delta, u^{\delta,*}
T^{\on{vert}}P(X))_{k,p,\delta} $.  Let
$$ \cH_{\delta} := \Omega^2(\Sigma,P(\g))_{{0,3}} \oplus
\Omega^0(\Sigma,P(\g))_{0,3} \oplus
\Omega^{0,1}(\Sigma^\delta, u^{\delta,*}
T^{\on{vert}}P(X))_{{0,3},\delta}$$
with norm
\begin{equation} \label{deltanorm} 
 \Vert (\phi,\psi,\eta) \Vert_{{\delta}}^2 = \Vert \phi
 \Vert_{{0,3}}^2 + \Vert \psi \Vert_{0,3}^2 + \Vert \eta
 \Vert_{0,3,\delta}^2 .\end{equation}
Let
$$ \cI_{\delta} := \Omega^1(\Sigma, P(\g))_{{1,3}} \oplus
\Omega^0(\Sigma^\delta, u^{\delta,*}
T^{\on{vert}}P(X))_{{1,3},\delta} $$
with norm
$$ \Vert (a, \xi ) \Vert_{{\delta}}^2 = \Vert a \Vert_{{1,3}}^2 + \Vert
\eta \Vert_{{1,3},\delta}^2 .$$
Locally the moduli space of polystable vortices is in bijection with
the zero set of the map
\begin{multline} \label{Fdel}
 \cF_{A,\mul{u}}^{{\delta}} : \cI_{{\delta}} \to \cH_{{\delta}} \\
(a,\xi) \mapsto \left( F_{A + a} + \Vol_\Sigma
 \exp_{\mul{u}^{{\delta}}}(\xi)^* P(\Phi), \d_{A,u_+}^* (a,\xi),
 \Psi_{\mul{u}^{{\delta}}}(\xi)^{-1} \ol{\partial}_{A + a}
 \exp_{\mul{u}^{{\delta}}}(\xi) \right) .\end{multline}
Here the second map enforces a slice condition.  That
$\cF_{A,\mul{u}}^{{\delta}}$ is well-defined follows from Sobolev
embedding: In particular, there is a $\delta$-uniform embedding
\begin{equation} \label{unifsob}
 \Omega^0(\Sigma^\delta, u^{\delta,*}
T^{\on{vert}}P(X))_{{1,3},\delta} \to \Omega^0(\Sigma^\delta,
u^{\delta,*} T^{\on{vert}}P(X))_{0,\infty}, \end{equation} 
since the dimensions of the cone in the cone condition \cite[Chapter
  4]{ad:so} are uniformly bounded in $\delta$, and the metric
uniformly comparable to the flat metric.

\begin{lemma} \label{approx} 
Let $(A,\mul{u})$ be a symplectic vortex on a nodal curve with a
single node $w = (w^+,w^-)$.  There exist constants $c_0,c_1 > 0$ such
that if
$|\delta| < c_1, \rho > 1/c_1 $ and$ |\delta| \rho^4 < c_1 $
then the pair $(A,\mul{u}^\delta) \in \A(P) \times \Gamma(P(X))$
satisfies
\begin{equation}
\Vert \cF_{A,\mul{u}}^{{\delta}}(0,0) \Vert \leq c_{0} |\delta|^{1/3}
\rho^{2/3}.
\end{equation}
\end{lemma}

\begin{proof}  
The expression $\olp_{A} \mul{u}^\delta$ can be expressed as a sum of
terms involving derivatives of the cutoff function $\alpha$, terms
involving derivatives of $\xi_j$, and terms involving the connection
$A$ on the bubble region.  The derivative of $\alpha$ is bounded by $
C/ \rho | \delta|^{1/2}$, while the norm of $\xi_j$ is bounded by
$C\rho |\delta |^{1/2}$ on the gluing region.  Hence the term
involving the derivative of $\alpha$ is bounded and supported on a
region of area less than $ C|\delta| \rho^2$.  In the given
trivialization we have
$$ \olp_A u =    \olp u + A_X^{0,1}(u) $$
where $A^{0,1}$ is the $0,1$-form defined by $A \in \Omega^1(B_R,\g)$
is the connection $1$-form in the local trivialization and
$A_X^{0,1}(u)$ is the corresponding form with values in $T^{\on{vert}}
P(X) \otimes_\R \C$.  We have
$$ \Vert A_X^{0,1}(u) \Vert_{0,3,\delta} \leq \Vert A_X^{0,1}(u)
 \Vert_{0,3}  $$
since $p \ge 2$; for $p = 2$ the $W^{0,3,\delta}$ and $W^{0,3}$ norms
are the same, by conformal invariance; for $p > 2$ the $0,3$-norm is
strictly greater.  Hence
$$ \Vert \olp_{A} u^\delta \Vert_{0,3,\delta} \leq C
\max(|\delta|^{1/3} \rho^{2/3},|\delta|). $$
The moment map term $F_{A} + (\mul{u}^\delta)^* P(\Phi) \Vol_\Sigma$
vanishes except on $|\kappa_+| \leq \rho |\delta^{1/2}|$, where it is
uniformly bounded.  Hence for $\delta$ small
$$ \Vert F_A + (\mul{u}^\delta)^* P(\Phi) \Vol_\Sigma \Vert_{0,3} \leq
C \rho^{2/3} |\delta|^{1/3} .$$
The statement of the lemma follows.  
\end{proof} 

We also wish to perform the gluing construction in families, that is,
for each nearby vortex and gluing parameter we wish to find a solution
to the vortex equations on $\Sigma^\delta$.  Define
\begin{multline} \label{Fdeldef}
 \cF_{A,\mul{u}}^{D,\delta} : \Def_\Gamma(A,\mul{u}) \times
 \cI_{{\delta}} \to \cH_{{\delta}}, \quad
(a,\mul{\xi},a_1,\xi_1) \mapsto \\ \left( F_{A} + a_0
   + a_1 + \Vol_\Sigma \exp_{\mul{u}^{{\delta}}}(\mul{\xi}_0^\delta +
 \xi_1)^* P(\Phi), \d_{A,\mul{u}^\delta}^* (a_0 +
 a_1,\mul{\xi}_0^\delta + \xi),
\right. \\ \left. \Psi_{\mul{u}^{{\delta}}}(\mul{\xi}_0^\delta + \xi_1)^{-1}
 \ol{\partial}_{A + a_0 + a_1}
 \exp_{\mul{u}^{{\delta}}}(\mul{\xi}_0^\delta + \xi_1)
 \right) \end{multline}

The following is proved in the same way as Lemma \ref{approx} and left
to the reader:  

\begin{lemma} \label{approx2} 
Let $(A,\mul{u})$ be as above.  There exist constants $c_0,c_1 > 0$
such that if $|\delta| < c_1, \rho > 1/c_1 $, $\Vert (a,\xi) \Vert
\leq c_1$ and $|\delta| \rho^4 < c_1 $ then
\begin{equation}
\Vert \cF_{A,\mul{u}}^{D,{\delta}}(a,\xi,0,0) \Vert \leq c_{0}
|\delta|^{1/3} \rho^{2/3}.
\end{equation}
\end{lemma}

\noindent {\em Step 2: Uniformly bounded right inverse}

In preparation for the construction of the uniformly bounded right
inverse of $\ti{D}_\delta$ we define the {\em intermediate family}
$(A,\mul{u}_0^\delta)$ of gauged holomorphic maps on the nodal curve
$\hat{\Sigma}$ is the family defined by the equations \eqref{approx},
using the identification of $\hat{\Sigma}$ and $\hat{\Sigma}^\delta$
away from the gluing region.  Thus $u_0^\delta$ is constant in a
neighborhood of the node $w^\pm$.  We identify $(\mul{u}_0^\delta)^*
T^{\on{vert}} P(X)$ with $\mul{u}^* T^{\on{vert}} P(X)$ by geodesic
parallel transport.

\begin{lemma} \label{linconverge}  
The operator $\ti{D}_{A,\mul{u}_0^\delta}$
converges in the operator norm to $\ti{D}_{A,\mul{u}}$ as $\delta \to
0, \rho \to \infty, \delta \rho^4 \to 0$.
\end{lemma}

\begin{proof}
The section $\mul{u}_0^\delta$ converges in the $W^{1,3}$ norm to $u$
as $\rho^2 |\delta|^{1/2} \to 0$.  It follows that the operator $\xi
\mapsto \Vol_\Sigma (u^\delta_0)^* L_\xi P(\Phi)$ converges to $\xi
\mapsto \Vol_\Sigma u^* L_\xi P(\Phi)$.  Hence $
\d_{A,u_0^\delta,\eps}$ converges to $\d_{A,u,\eps}$, and similarly
for $\d_{A,u_0,\eps}^*$.  The operator $ D_{A,\mul{u}_0^\delta}$
converges to $D_{A,\mul{u}}$, as in Lemma \ref{rightex}.
\end{proof}

\begin{proposition}
\label{rightvor}
Let $(A,\mul{u})$ be a nodal vortex, and $(A,\mul{u}^\delta)$ the
approximate solution constructed above.  There exist constants
$c,C > 0$ such that if $| \delta | < c$ then there exists an
approximate right inverse $T_{\delta}$ of the parametrized linear
operator $\ti{D}_{\delta} := \ti{D}_{A,\mul{u}^\delta}$ that is,
a map $T_{\delta} : \cH_{\delta} \to \cI_{\delta} $
such that
$$ \Vert (\ti{D}_{\delta} T_{\delta} - I) \eta \Vert_{\delta} \leq \hh
\Vert \eta \Vert_{\delta}, \quad \Vert T_{{\delta}} \eta
\Vert_{\delta} \leq C \Vert \eta \Vert_{\delta} .$$
\end{proposition} 

Given such an approximate inverse, we obtain a uniformly bounded right
inverse $Q_{\delta}$ to $\ti{D}_{\delta}$ by the formula
$$ Q_{\delta} = T_{\delta} (\ti{D}_{\delta} T_{\delta})^{-1} = \sum_{k \ge 0}
T_{\delta} (\ti{D}_{\delta} T_{\delta} - I )^k .$$

\begin{proof}[Proof of \ref{rightvor}] 
By the regularity assumption, $\ti{D}_{A,\mul{u}} $ is surjective when
restricted to the space of vectors $(a,\mul{\xi})$ such that $\xi_0(0)
= \xi_1(\infty)$.  By Lemma \ref{linconverge}, $\ti{D}^0_\delta :=
\ti{D}_{A,\mul{u}_0^\delta}$ is surjective for sufficiently small
$\rho,\delta$, when restricted to the same space.  The approximate
right inverse is constructed by composing a cutoff operator
$K_\delta$, right inverse $Q_\delta$, and gluing operator $R_\delta$,
as follows.  (In other words, the right inverse is defined by
truncating the given functions, applying the right inverse for the
linearized operator on the nodal curve, and then gluing back together
using cutoff functions again.)  For simplicity we assume that there is
a single node.  Define the {\em cutoff operator}
\begin{multline}
 K_{\delta}: 
\Omega^{0,1}(\hat{\Sigma}^\delta,\mul{u}^{\delta,*} T^\vert
P(X))_{0,3,\delta} 
\to 
 \Omega^{0,1}(\hat{\Sigma},(\mul{u}^\delta_0)^{*}
T^{\on{vert}}P(X))_{0,3} \\
K_\delta(\eta) =
\begin{cases}  0  & \kappa_\pm( z) \in B_{| \delta|^{\hh}} (0)  \\
                \eta & \text{otherwise}
                \end{cases} .\end{multline}
Then $ \Vert K_\delta(\eta) \Vert_{0,3} \leq \Vert \eta
\Vert_{0,3,\delta} .$ Define the {\em gluing operator}
$$ R_{\delta}: 
\Omega^0(\hat{\Sigma},(\mul{u}^\delta_0)^{*} TX)_{{1,3}} \to
\Omega^0(\hat{\Sigma}^{\delta},u^{\delta,*}
T^{\on{vert}}P(X))_{{1,3},\delta}
$$
as follows.  Let $\hat{\Sigma}_\pm^*$ denote the complements of small
balls around the nodes $ \hat{\Sigma}_\pm^* = \Sigma_\pm - B_{\rho^2 |
  \delta|^{1/2}} (w^\pm)$.  Let $\pi_\pm : \hat{\Sigma}_\pm^* \to
\hat{\Sigma}^\delta $ denote the inclusions.  These induce maps of
sections with compact support in $\hat{\Sigma}_\pm^*$, $ \pi_{\pm,*}:
\, \Omega^0_c( \hat{\Sigma}_\pm^*,u_\pm^* TX) \to
\Omega^0(\hat{\Sigma}^\delta, u^{\delta,*} TX ) .$ Define $
R_\delta(\xi) = \xi^\delta $ where
$$ \xi^\delta = \pi_{+,*} \beta_{\rho,\delta} (\xi_+ - \xi_+(w_+)) +
\pi_{-,*} \beta_{\rho,\delta} (\xi_- - \xi_-(w_-)) + \xi(w) $$
for $ \kappa_\pm(z) \in B_{| \delta|^{1/2} \rho^2} (0) $ and $
\xi^\delta = \xi$ otherwise.  Here $x = \xi^\pm(w^\pm)$ is the value
of $\xi$ at the node.  Define $ T_\delta := (I \times R_\delta)
Q_\delta (I \times K_\delta) .$ That is, if $ (a,\xi) = Q_\delta
K_\delta(\phi,\psi,\eta) $ then $ T_\delta(\phi,\psi,\eta) =
(a,\xi^\delta) .$ The map $T_\delta$ is the required approximate right
inverse.  The difference $(\ti{D}_\delta T_\delta -
I)(\phi,\psi,\eta)$ is the sum of terms
$$ \d_{A,u^\delta} (a,\xi^\delta) - \phi, \quad
\d_{A,u_0}^* (a,\xi^\delta) - \psi, \quad
D_{A,u^\delta} (a,\xi^\delta) - \eta .$$
By definition $\d_{A,u_0^\delta}(a,\xi) = \phi$, so the first
difference has contributions involving the difference between
$u_0^\delta$ and $u^\delta$, and between $\xi^\delta$ and $\xi$.
But since $\d_{A,u}(a,\xi) = \d_A a + u^* L_{\xi_X} P(\Phi)$, these
terms are zeroth order in $u,\xi$, 
$$  \Vert (u^\delta)^* L_{\xi^\delta_X} P(\Phi) - u_0^\delta L_{\xi_X}
P(\Phi) \Vert_{0,3} \leq  C | \delta | \Vert \xi \Vert_{C^0} ;$$
that is, a constant times the area of $\pi_{-,*} (\hat{\Sigma}_-^*)$,
which goes to zero as $\delta$ does.  A similar discussion holds for
the second difference.  The third difference has terms arising from
the cutoff function and the term $a_X^{0,1}(u)$ on the bubble region
$|\kappa_+(z)| \leq | \delta|^{1/2}$.  
We have
$$ \Vert a_X^{0,1}(u) \Vert_{0,3,\delta} \leq 
\Vert a_X^{0,1}(u) \Vert_{0,3} \leq C |\delta|  .$$
since $p \ge 2$; for $p = 2$ the $W^{0,3,\delta}$ and $W^{0,3}$ norms
are the same, by conformal invariance; for $p > 2$ the $0,3$-norm is
strictly greater.  The term involving the derivative of the cutoff
function satisfies
$$ \Vert \d \beta_{\rho,\delta} (\xi - \xi(w)) \Vert_{0,3,\delta} \leq 
c \log( \rho^2 )^{-2/3} \Vert (\xi - \xi(w)) \Vert_{1,3} $$
by Lemma \ref{careful}, and so vanishes in the limit $\rho \to
\infty$.  Using the uniform bound on $Q_\delta$, the total difference
is bounded by $C ( \log( \rho^{2})^{-2/3} + |\delta|) \Vert (\phi,
\psi,\eta) \Vert $, and so vanishes in the limit $\delta \to 0,\rho
\to \infty, |\delta| \rho^2 \to 0$.  The uniform bound on $T_{\delta}$
follows from the uniform bound on $Q_\delta$ and the cutoff and
extension operators.
\end{proof}

\noindent {\em Step 3: Uniform quadratic estimate}

\begin{proposition} \label{quad}
There exist constants $c,C > 0$ such that if $\Vert \xi
\Vert_{1,3,\delta} < c$, $ |\delta| < c, \rho > 1/c$ and $ |\delta|
\rho^4 < c $ then the map $\cF^\delta_{A,\mul{u}}$ satisfies a
quadratic bound
$$ \Vert D \cF_{A,\mul{u}}^\delta(a_1,\xi_1) -
\ti{D}_{A,\mul{u}^\delta}(a_1,\xi_1) \Vert_{\delta} \leq C \Vert a,\xi
\Vert_{\delta} \Vert a_1, \xi_1 \Vert_{\delta} .$$
\end{proposition} 

\begin{proof} 
The norm of the non-linear part of the curvature $\Vert [a,a_1]
\Vert_{{0,3}}$ is bounded by Sobolev multiplication.  The other term
appearing in the first vortex equation satisfies
$$ \Vert  \exp_{\mul{u}^\delta} (\xi_0)^*  L_{\xi_{1,X}}
P(\Phi) - (\mul{u}^\delta)^* L_{\xi_{1,X}}
P(\Phi) \Vert_{{0,3}} \leq C \Vert \xi_0 \Vert_{1,3,\delta} \Vert
\xi_1 \Vert_{1,3,\delta} $$
for some constant $C$ independent of $\delta$, using that
$W^{1,3,\delta}$ norm controls the $W^{0,3}$ norm uniformly.  The
non-linear terms in the Cauchy-Riemann equation are estimated as in
Theorem \ref{quad} and \cite[Section 3.5, Lemma 10.3.1]{ms:jh}; note
that we are fixing the complex structure on $\Sigma$, which avoids the
more complicated analysis we gave in the previous section.  The second
vortex equation also involves a term of mixed type $ \Psi_{u}(\xi +
\xi_1)^{-1} (a_1)_X^{0,1}(\exp_u(\xi + \xi_1)) - \Psi_u(\xi)^{-1}
(a_1)_X^{0,1}(\exp_u(\xi)) .$ It follows from uniform Sobolev
embedding that this difference has $0,3,\delta$-norm bounded by $C
\Vert a_1 \Vert_{{1,3}} \Vert \xi_1 \Vert_{1,3,\delta}$ for some
constant $C$ independent of $\delta$.
\end{proof} 

\noindent {\em Step 4: Implicit Function Theorem}

\begin{theorem}
\label{gluevor} Let $(A,\mul{u})$ be a regular stable 
 nodal vortex of combinatorial type $\Gamma$.  There exist constants
 $\eps_0,\eps_1 > 0$ such that for every $(a,\mul{\xi}, \delta) \in
 \Def_\Gamma(A,\mul{u})$ with norm less than $\eps_0$, there exists
 a unique $(\phi,\psi,\eta)$ of norm less than $\eps_1$ such that if
 $(a_1,\xi_1) = Q_\delta (\phi,\psi,\eta)$ then $(A + a_0 + a_1,
 \exp_{\mul{u}^{\delta}}(\mul{\xi}_0^{\delta} + \xi_1))$ is a
 symplectic vortex in Coulomb gauge with respect to
 $(A,\mul{u}^\delta)$.  The solution depends smoothly on
 $a_0,\mul{\xi}_0$, and transforms equivariantly the action of
 $\G(P)_{A,\mul{u}}$ on $\Def_\Gamma(A,\mul{u})$.
\end{theorem}  

\begin{proof} 
Uniform error and quadratic estimates are those for
$\cF_{A,\mul{u}}^\delta$ in Lemmas \ref{approx}, \ref{rightvor}, and
\ref{quad}, in a uniformly bounded neighborhood of $0$ in
$\Def_\Gamma(A,\mul{u})$.  Then the first claim is an application of
the quantitative version of the implicit function theorem (see for
example \cite[Appendix A.3]{ms:jh}).  Equivariance follows from
uniqueness of the solution given by the implicit function theorem,
since the map $\cF^{D,\delta}_{A,\mul{u}}$ is equivariant for the
action of $\G(P)_{A,\mul{u}}$.
\end{proof}

\noindent {\em Step 5: Rigidification}

As in the case of holomorphic maps in the previous section, there is a
more natural way of parametrizing nearby symplectic vortices which
involves examining the intersections of the sections with submanifolds
of $P(X)$, and framings induced by parallel transport.  First we study
the differentiability of the evaluation maps.  The gluing construction
of the previous step gives rise to a deformation $(A_S,\mul{u}_S)$ of
$(A,\mul{u})$ with parameter space a neighborhood $S$ of $0$ in
$\Def(A,\mul{u})$, and so a map
$ S \to \ol{M}_n(\Sigma,X), \ s \mapsto (\hat{\Sigma}_s,A_s,\mul{u}_s) $
Consider the map
\begin{equation} \label{evs}
\ev: (\hat{\Sigma} - U) \times S \to P(X), \quad (z,s) \mapsto
\mul{u}_s(z) .
\end{equation} 

\begin{proposition} \label{diffv} 
The map $\ev$ of \eqref{evs} is $C^1$ for the family constructed by
gluing in Theorem \ref{gluevor} using the exponential gluing profile.
\end{proposition} 

\begin{proof}  We denote by $\mul{u}_S^\pre: \hat{\Sigma}_S \to X$ 
the family obtained by pre-gluing only, that is, omitting the step
which solves for an exact solution.  We denote by $\ev^\pre$ the map
$$ \ev^\pre : (\hat{\Sigma} - U) \times S \to P(X), \quad (z,s) \mapsto
\mul{u}_s^\pre(z) .$$
This map is independent of the gluing parameters, and is therefore
$C^1$.  We write $s = (a_0,\xi_0)$ and $A_s = A + a_0 + a_1, \mul{u}_s
= \exp_{\mul{u}^{\pre}_s}(\xi_0^\delta + \xi_1)$. The corrections
$a_1,\xi_1$ depend smoothly on $a_0,\xi_0$, by the implicit function
theorem, and so $\xi_1(z)$ depends smoothly on $a_0,\xi_0$.  Next we
take the derivative with respect to the gluing parameter. Let
$(A,\mul{u})$ be a nodal symplectic vortex, $(A,\mul{u}^\delta)$ the
pre-glued pair (we omit the parameter $\rho$ controlling the diameter
of the gluing region from the notation) and consider the equation $
\cF_{A,\mul{u}^\delta} (a_0 + a_1, \xi_0^\delta + \xi_1) = 0 .$
Let $\ti{D}_\delta$ denote the derivative of $\cF_{A,\mul{u}^\delta}$.
Differentiating with respect to $\delta$ gives
$$ \ti{D}_\delta \left(\ddd a_1 , D\exp_{\mul{u}^\delta}(\xi_0^\delta;
0 , \ddd \xi_1 )\right) = - \ti{D}_\delta \left(0,
D\exp_{\mul{u}^\delta}(\xi_0^\delta; \ddd \mul{u}^\delta, \ddd
\xi_0^\delta)\right) .$$
The same arguments as in the proof of Theorem \ref{diffev} show that
there exists a constant $C > 0 $ such that the right hand side is
bounded in norm by $C e^{-1/\delta}$.  On the other hand, the norm of
the left-hand side $\ti{D}_\delta$ is uniformly bounded from below in
terms of the norm of $\ddd a_1, \ddd \xi_1$, by the quadratic
estimates.  It follows that $(\ddd a_1, \ddd \xi_1)$ is also bounded
in norm by $C e^{-1/\delta}$.  Hence $ \lim_{\delta \to 0}
\partial_\delta \ev = 0 $.  It follows that $D\ev$ has a continuous
limit as ${\delta} \to 0$.
\end{proof}  

Choose a path $\gamma: [0,1] \to \Sigma$ in the principal component
and an element $\phi_0 \in P_{\gamma(0)}$.  Let $\tau_\gamma(A):
P_{\gamma(0)} \to P_{\gamma(1)}$ denote parallel transport.  By an
{\em $m$-framed} family of marked curves, we mean a family of curves
together with an $m$-tuple of points in $P$.  Given a family
$(\hat{\Sigma}_S,A_S,\mul{u}_S)$ of gauged holomorphic maps over a
parameter space $S$, a collection of codimension two submanifolds
$\mul{Y} = (Y_1,\ldots,Y_k)$ in $P(X)$, and a collection of paths
$\mul{\gamma} = (\gamma_1,\ldots, \gamma_l)$ with the same initial
point $y_0$ to $y_j, j = 1,\ldots, l$, define a family of marked,
framed curves
$\hat{\Sigma}_S^{\mul{Y},\mul{u},\mul{\gamma},A} \to S$ 
by requiring that the additional marked points $z_{n+i}$ map to $Y_i$,
and the framings are given by parallel transport along the paths
$\gamma_i$.

\begin{definition} \label{compat2}
The data $\mul{Y},\mul{\gamma},A,\mul{u}$ are {\em compatible} if
\begin{enumerate}
\item each $Y_j$ intersects $u_j$ transversally in a single point $z_j
\in \hat{\Sigma}$;
\item if $(a,\xi) \in \ker \ti{D}_{A,\mul{u}}$ satisfies $\xi(z_{n+j})
  \in T_{\mul{u}(z_{n+j})} P(X) $ for $j = 1,\ldots, k$ and $D_A
  \tau_{\gamma_i}(a) = 0 $ for $i = 1,\ldots, l$ then $(a,\xi) = 0$.
\item the curve $\hat{\Sigma}$ marked with the additional points
  $z_{n+1},\ldots, z_{n+k}$ is stable.
\item if some automorphism of $(\hat{\Sigma},\mul{u})$ maps $z_i$ to
$z_j$ then $Y_i$ is equal to $Y_j$.
\end{enumerate} 
\end{definition} 

The second condition says that there are no infinitesimal deformations
which do not change the positions of the extra markings or framings.  

\begin{proposition} \label{versalrigid2}
  Let $(A,\mul{u})$ be a parametrized regular stable nodal vortex, and
  $(A_S, \mul{u}_S) \to S$ the stratified-smooth universal deformation
  constructed by the gluing construction.  There exists a collection
  $(\mul{Y},\mul{\gamma})$ compatible with $(\mul{u},A)$. Furthermore, if
  $(\mul{Y},\mul{\gamma})$ is compatible with $(A,\mul{u})$, then the
  rigidified family $\hat{\Sigma}_S^{\mul{Y},\mul{u},\mul{\gamma},A} \to
  S$ of \eqref{rigidify} is a stratified-smooth deformation of the
  marked-curve-with-framings
  $\hat{\Sigma}^{\mul{Y},\mul{u},\mul{\gamma},A}$ which defines an
  immersion of $S$ into the parameter space for the universal
  deformation of the central fiber.
\end{proposition}

\begin{proof} 
First we show the existence of a compatible collection.  Suppose that
the second condition is not satisfied for some $(a,\xi)$.  Suppose
first that $\xi \neq 0$.  Let $z_{n+1}$ be a point with $\xi(z_{n+1})
\neq 0$, and choose a codimension two submanifold $Y_{n+1}$ transverse
to $u$ near $u(z_{n+1})$, and such that $TY_{n+1}$ does not contain
$\xi(z_{n+1})$.  Adding $Y_{n+1}$ to the list of submanifolds
decreases the dimension of the space of $(a,\xi)$ satisfying the
condition in (b) by at least one.  Repeating this process, we may
assume that the only elements satisfying the condition in (b) have
$\xi =0 $.  Suppose that $\xi$ is zero, so that $a$ is necessarily
non-zero.  Choose an additional marked point $y_{l+1}$ and a path
$\gamma_{l+1}$ from the base point $y_0$ to $y_{l+1}$ such that the
derivative of the parallel transport over $\gamma$ with respect to $a$
over is non-zero.  Appending $\gamma_{l+1}$ to the list of path
decreases the dimension of $(a,\xi)$ satisfying the condition in (b)
by at least one.  Hence the process stops after finitely many steps,
after which the kernel is trivial.  The proof of the second claim is
similar to Proposition \ref{versalrigid} and will be omitted.
\end{proof} 

\noindent {\em Step 6: Surjectivity}

We show that any nearby vortex appears in the family constructed
above.  First, we show:

\begin{proposition} \label{vorsurj} 
Let $(A,\mul{u})$ be a regular strongly stable symplectic vortex.
There exists a constant $\eps > 0$ such that if $(A_1,u_1) = (A + a,
\exp_{\mul{u}^\delta}(\xi))$ with $\Vert a \Vert_{1,3} + \Vert \xi
\Vert_{1,3,\delta} \leq \eps$ then after gauge transformation we have
$(A_1,u_1) = (A+ a_0 + a_1,\exp_{\mul{u}^\delta}(\xi_0^\delta +
\xi_1))$ for some $(a_0,\xi_0) \in \ker(\ti{D}_{A,\mul{u}})$ and
$(a_1,\xi_1)$ in the image of $Q_\delta$.
\end{proposition}  

\begin{proof}
We claim that for some constant $C > 0$, we have $(a,\xi) =
(a_0,\xi_0^{\mul{\delta}}) + (a_1,\xi_1)$ for some $(a_0,\xi_0) \in
\ker \ti{D}_{A,\mul{u}}$ and $(a_1,\xi_1) \in \Im
\ti{D}_{A,\mul{u}^{{\delta}}}^*$ with norm $ \Vert (a_1,\xi_1) \Vert
\leq C \Vert (a_1,\xi_1) \Vert.$ Given the claim, the proposition
follows by the uniqueness statement of the implicit function theorem.
For any $c > 0$ there exists $\delta_0$ such that for $\delta <
\delta_0$,
$$ \Vert \ti{D}_{A,\mul{u}^{{\delta}}} (
a_0^{{\delta}},\mul{\xi}_0^{{\delta}}) \Vert \leq c \Vert
(a_0^{{\delta}},\mul{\xi}_0^{{\delta}}) \Vert $$
by estimates similar to those of Lemma \ref{approx}.  Thus the image
of $\ker \ti{D}_{A,\mul{u}}$ is transverse to $\Im
\ti{D}_{A,\mul{u}^\delta}^*$, for $\delta$ sufficiently small, since
it meets $\Im \ti{D}_{A,\mul{u}^\delta}^*$ trivially and projects
isomorphically onto $\ker \ti{D}_{A,\mul{u}^\delta}$, by gluing for
indices, as in \cite[Theorem 2.4.1]{orient}.  The claim then follows
from the inverse function theorem.
\end{proof} 

Given a strongly stable symplectic vortex $(A,\mul{u})$ with stable
domain $\hat{\Sigma}$, let $(A_S,\mul{u}_S)$ be the family given by
the gluing construction above.  Otherwise, if $\mul{Y}$ is not stable,
let $\mul{Y} = (Y_1,\ldots, Y_l)$ be a collection of codimension two
submanifolds of $P(X)$, and consider the family
$(A,\mul{u}^{\mul{Y}})$ with additional marked points given by
requiring that the additional marked points $z_{n+i}$ map to $Y_i$.
Let $(A_S,\mul{u}_S)$ denote the family obtained by applying the
gluing construction for $(A_S,\mul{u}_S^{\mul{Y}})$, and then
forgetting the additional marked points.

\begin{lemma}   Suppose that $(A_i,u_i)$ Gromov converges to $(A,u)$. 
After a sequence of gauge transformations, for any $\eps$, there
exists $i_0$ such that if $i > i_0$ then there exists
$\delta,(a_i,\xi_i)$ satisfying $(A_i,u_i) = (A + a_i,
\exp_{\mul{u}^\delta}(\xi_i))$ with $\Vert a_i \Vert_{1,3} + \Vert
\xi_i \Vert_{1,3,\delta} \leq \eps$.
\end{lemma}

\begin{proof}  By definition of Gromov convergence, after gauge transformation  $A_i$ $C^0$-converges to $A$
and converges uniformly in all derivatives on the complement of the
bubbling set \cite{ott:thesis}.  The exponential decay estimate
\cite[Lemma A.2.2]{ott:thesis} show that $u_i$ converges to $u$ on the
complement of the nodes, uniformly in all derivatives on compact sets,
and whose derivative on the gluing region is uniformly bounded in the
$\delta$-dependent metric. It follows that $u_i =
\exp_{u^\delta}(\xi_i)$ for some $\delta$ and $\xi_i \in
\Omega^0(\Sigma^\delta, (u^\delta) T^{\on{vert}} P(X))$ with $\Vert
\xi_i \Vert_{1,3,\delta} < \eps$.  To obtain the improved convergence
for the connection, note that $ F_{A_i} + (u_i)^* P(\Phi) = 0 $ and
the corresponding equations for the limit $(A,u)$ imply that
$$ F_{A_i} - F_A = \d_A (A_i - A) - (A - A_i) \wedge (A - A_i) =
(u_i)^* P(\Phi) - u^* P(\Phi) .$$
Since $u_i^* P(\Phi)$ is bounded and converges to $u^* P(\Phi)$ on the
complement of the bubbling set, and $A_i$ converges to $A$ in $C^0$
hence $W^{0,3}$, the right hand side converges to $0$ in $W^{0,3}$ as
$i \to \infty .$ After gauge transformation we may assume that $
\d_A^* (A - A_i) = 0$.  Then the elliptic estimate for the operator
$\d_A \oplus \d_A^*$ implies that $A - A_i$ converges to zero in
$W^{1,3}$.
\end{proof} 

\begin{corollary}  \label{versalvor}
$(A_S,\mul{u}_S)$ is a stratified-smooth versal deformation of
  $(\hat{\Sigma},A,\mul{u})$.
\end{corollary}

\begin{proof} 
Proposition \ref{vorsurj} implies that any family
$(A_{S^1}^1,\mul{u}_{S^1}^1)$ is obtained by pull-back from
$(A_S,\mul{u}_S)$, in case $\hat{\Sigma}$ is stable, or obtained from
the family obtained by adding the marked points mapping to
submanifolds, in general.
\end{proof} 

\noindent {\em Step 7: Injectivity}

We show that any nearby vortex appears once in our family, up to the
action of $\Aut(A,u)$; this is part of the following:

\begin{theorem} Any family $(A_S,\mul{u}_S)$ constructed
by gluing using the exponential gluing profile is a strongly universal
stratified-smooth deformation of $(A,\mul{u})$.
\end{theorem} 

\begin{proof}  
Let $\ol{Z}_n(P,X)$ denote the moduli space of marked symplectic
vortices up to equivalences that involve only the identity gauge
transformation, so that $ \ol{M}_n(P,X) = \ol{Z}_n(P,X) / \G(P) .$ Let
$(A,\mul{u})$ be a stable marked vortex, and $W_{A,\mul{u}}$ a slice
for the gauge group action on $\ol{Z}(P,X)$, so that
$$W_{A,\mul{u}}/\G(P)_{A,\mul{u}} \to \ol{M}_n(P,X) $$
is a homeomorphism onto its image.  Let $\Aut_0(A,\mul{u})$ denote the
subgroup of $\Aut(A,u)$ acting trivially on $P$, so that
$\G(P)_{A,\mul{u}} = \Aut(A,u)/\Aut_0(A,\mul{u})$ is the stabilizer of
$(A,\mul{u})$ under the gauge action.  Let $(A_S,\mul{u}_S)$ denote a
universal deformation of $(A,\mul{u})$ constructed by gluing using the
exponential gluing profile.  We claim that the map
\begin{equation} \label{glueS3} S/\Aut_0(A,\mul{u}) \to W_{A,\mul{u}}, \quad [s]
\mapsto [A_s,\mul{u}_s] \end{equation}
is an injection.  Indeed, rigidification produces an injection
\begin{equation} \label{glueS4} S/\Aut_0(A,\mul{u}) 
\to \ol{M}_{n+k,l}(\Sigma)/\Aut_0(A,\mul{u}), \quad [s] \mapsto
[\Sigma^{A_s,\mul{u}_s,\mul{Y},\mul{\gamma}}] \end{equation}
where $\Aut_0(A,\mul{u})$ acts by re-ordering the marked points.
Since this map factors through \eqref{glueS3}, the claim follows.  If
$(A_{S^1},\mul{u}_{S^1}^1)$ is a family of symplectic vortices giving
a deformation of any fiber of $(A_S,\mul{u}_S)$< then Corollary
\ref{versalvor} together with injectivity shows that this family is
obtained by pull-back by some map $S^1 \to S$.  Hence
$(A_S,\mul{u}_S)$ is a stratified-smooth strongly universal
deformation of $(A,\mul{u})$.
\end{proof} 

\begin{theorem}  \label{regvortex}  Let $X$ be a Hamiltonian $G$-manifold
equipped with a compatible invariant almost complex structure $J \in
\J(X)^G$.  The maps
\begin{equation} S \to \ol{M}_n(\Sigma,X), \quad s
\mapsto [A_s,\mul{u}_s] \end{equation}
associated to the universal deformations constructed above equip the
locus $\ol{M}^\reg_n(\Sigma,X)$ of regular stable symplectic vortices
with the structure of a stratified-smooth orbifold.  If the local
coordinates near the nodes are chosen compatibly and the gluing
profile is the exponential gluing profile, then the deformations
provide $\ol{M}^\reg_n(\Sigma,X)$ with the structure of a
$C^1$-orbifold.
\end{theorem}  

\begin{proof}   
It suffices to show that the charts given by two sets
$\mul{Y}_j,\mul{\gamma}_j$ are compatible.  Define $\mul{Y} =
\mul{Y}_1 \cup \mul{Y}_2$ and $m = m_1 + m_2$ the total number of
extra points.  Similarly let $\mul{\gamma}$ be the union of
$\mul{\gamma}_1$ and $\mul{\gamma_2}$ of total number $l = l_1 +
l_2$. The family $ \hat{\Sigma}^{\mul{Y},\mul{u},\mul{\gamma},A}_S$
admits a proper \'etale forgetful map $
\hat{\Sigma}_S^{\mul{Y},\mul{u},\mul{\gamma},A} \to
\hat{\Sigma}_S^{\mul{Y}_j,\mul{u},\mul{\gamma}_{j},A}, \quad j = 1,2$
whose fiber consists of the re-orderings of the points for $\mul{Y}$
induced by automorphisms of $\Aut(A,\mul{u})$ that fix the ordering
for $\mul{Y}_j$.  It follows that the corresponding charts are
$C^1$-compatible.
\end{proof}  

\begin{remark} As discussed in Remark \ref{smooth}, the Theorem 
implies that if $\ol{M}^{\reg}_n(\Sigma,X)$ is compact then it admits
a (non-canonical) smooth structure.
\end{remark} 

 Let $\ol{M}_n(\Sigma)$ denote the moduli space of stable maps to
 $\Sigma$ with homology class $[\Sigma]$, $n$ markings and genus that
 of $\Sigma$, or in other words, {\em parametrized} stable curves with
 principal component isomorphic to $\Sigma$.  Forgetting the pair
 $(A,u)$ gives a forgetful morphism $ \ol{M}^{\reg}_n(\Sigma,X) \to
 \ol{M}_n(\Sigma) .$ Using the differentiable structure defined above,
 the evaluation maps are differentiable but unfortunately the
 forgetful morphisms are not, unless one uses a different gluing
 profile for the moduli space of vortices with one less marking.  More
 precisely, the forgetful morphism $\ol{M}^\reg_n(\Sigma,X) \to
 \ol{M}_n(\Sigma)$ is continuous and $C^1$ near any pair $(A,\mul{u})$
 whose domain is stable as an element of $\ol{M}_n(\Sigma)$, and a
 submersion near the boundary of $\ol{M}_n(\Sigma)$.  For the standard
 smooth structure on $\ol{M}_n(\Sigma)$, the forgetful morphism
 $\ol{M}^\reg_n(\Sigma,X) \to \ol{M}_n(\Sigma)$ is smooth.

The gluing construction has various parametrized versions.  For
example, in \cite{cross} we consider a moduli space of {\em polystable
  polarized vortices}, which consist of a vortex together with a lift
of the connection to the Chern-Simons line bundle.  In each of these
cases one applies the implicit function theorem using the linearized
operator for the parametrized problem to prove that any {\em
  parametrized regular} polystable vortex has a strongly universal
deformation in the parametrized sense.  In particular, any regular
polystable polarized vortex has a strongly universal deformation etc.

\def\cprime{$'$} \def\cprime{$'$} \def\cprime{$'$} \def\cprime{$'$}
  \def\cprime{$'$} \def\cprime{$'$}
  \def\polhk#1{\setbox0=\hbox{#1}{\ooalign{\hidewidth
  \lower1.5ex\hbox{`}\hidewidth\crcr\unhbox0}}} \def\cprime{$'$}
  \def\cprime{$'$}

\end{document}